%% file: parahoric_arxiv_v1.tex
\documentclass[10pt,a4paper, twoside, notitlepage]{article}
\usepackage{import} 
\import{}{preamble}

\bibliography{bibli.bib} 

\usepackage{libertinus}

\usepackage{fancyhdr}
\title{\uppercase{\textbf{\large{Reductification of parahoric group schemes}}}}
\author{\uppercase{Arnab Kundu}}
\date{}
\usepackage{doi}
\begin{document}
\maketitle
\abstract{Parahoric group schemes are certain possibly non-reductive, smooth, affine integral models of reductive group schemes defined over a henselian discretely valued field $K$ whose residue field is perfect. We show that any such group scheme $\mathscr{P}$ becomes reductive, in a particular regard, after a (possibly wildly ramified) finite Galois extension $L/K$. More precisely, we prove that there exists a reductive integral model $\mathscr{G}$ of the base change $\mathscr{P}_L$ such that $\mathscr{P}$ can be recovered as the smoothening of the subgroup of Galois invariants of the Weil restriction of $\mathscr{G}$. Our work extends results of Balaji--Seshadri and Pappas--Rapoport from the tamely ramified and simply-connected semisimple setting. 
	\par As an application, we establish a parahoric analogue of the Grothendieck--Serre conjecture in sufficiently good residue characteristics. Specifically, we confirm that generically trivial parahoric torsors are trivial whenever the generic reductive group is simply-connected. The proof proceeds by reducing the problem to a statement about a stacky reductive group over a stacky discrete valuation ring.}
\tableofcontents
\section{Analogue of the Grothendieck--Serre conjecture for parahoric groups.}
In the theory of torsors under reductive group schemes, the conjecture of Grothendieck--Serre enjoys a central place due to the simplicity of its statement while simultaneously having far reaching consequences (see \cite{ces_torsors} for a recent survey). It has a rich history and many important cases of this conjecture are currently known to be true; furthermore, no counterexample has been found yet. Given its success, it is natural to wonder if an extension of this conjecture might hold beyond the class of reductive group schemes. 
\par In this light, we begin our discussion with the following question that appeared in \ccite{bayer-fluckiger_first_parahoric}{Question~6.4}. A positive answer would provide an obvious extension of the Grothendieck--Serre-style results to parahoric group schemes in the sense of Definition~\ref{defn:parahoric}, since these group schemes generalise the reductive ones.
\par Before stating the question, we emphasise that throughout this article, all parahoric group schemes--including the reductive ones--are assumed to be fibrewise connected by default.
\begin{quest}[Bayer-Fluckiger--First]\label{quest:bayer_fluckigar_first}
	Given a parahoric group scheme $\scr{P}$ in the sense of Definition~\ref{defn:parahoric}\eqref{defn:parahoric:2} over a semilocal principal ideal domain $R$ with perfect residue fields and fraction field $F$, do we have an equality \bud \ker(H^1(R,\scr{P})\xrightarrow{\rho} H^1(K,\scr{P}))=\{\ast\}~?\eud
\end{quest}	
We note, however, a subtle difference between Question~\ref{quest:bayer_fluckigar_first} and \ccite{bayer-fluckiger_first_parahoric}{Question~6.4}. In our formulation, we inquire if the restriction $\rho$ to the generic point has trivial kernel, as opposed to \citeauthor{bayer-fluckiger_first_parahoric}, who asked if $\rho$ is actually injective. Although, by translating the origin of $H^1(R,\scr{P})$ (see, for example, the argument in \ccite{phd-thesis}{Remark~2.0.2}), this stronger condition of injectivity follows provided we have an affirmative answer to Question~\ref{quest:bayer_fluckigar_first} for each parahoric $\scr{P}'$ whose generic fibre is an inner form of the generic fibre of $\scr{P}$. 
\subsubsection*{Known answers to Question~\ref{quest:bayer_fluckigar_first}}
\bun 
\item The case when $\scr{P}$ is reductive over $R$ is a part of the Grothendieck--Serre conjecture for torsors under reductive groups. In this case, the question has an affirmative answer even when $R$ is a semilocal Dedekind domain, thanks to the work of Nisnevich in his thesis \cite{nisnevich_thesis} and of Guo in his Master's thesis \cite{ning_dedekind}. Consequently, it suffices to treat the case of non-reductive $\scr{P}$.
\item When $R$ is a complete discrete valuation ring, \citeauthor{bruhat_tits_buildings_III} had already established a partial positive result in \ccite{bruhat_tits_buildings_III}{lemme~3.9}. Moreover, in the case where the generic fibre of $\scr{P}$ is the group of automorphisms of a hermitian form, \citeauthor{bayer-fluckiger_first_parahoric} themselves asserted in \ccite{bayer-fluckiger_first_parahoric}{\textsection6} that Question~\ref{quest:bayer_fluckigar_first} has an affirmative answer. 
\item More recently, \citeauthor{zidani_gs} in \cite{zidani_gs, zidani_gs_parahoric_henselian_dvr_case} has provided a positive response when the generic fibre of $\scr{P}$ is a semisimple group scheme that is either simply-connected, or quasi-split and adjoint, and $R$ is even a semilocal Dedekind ring. 
\eun
Our aim in this article is to provide a new demonstration of the following result, addressing the case where $R$ is a henselian discrete valuation ring.

\but{A}[Corollary~\ref{cor:general_case_gs}]\label{thm:A}
Given a henselian discrete valuation ring $\ca{O}_K$ with a perfect residue field of characteristic $p\neq 2,3,5$ and field of fractions $K$, for a parahoric group $\ca{O}_K$-scheme $\scr{P}$ whose generic fibre $\scr{P}_K$ is a simply-connected semisimple group scheme such that $p\nmid (n+1)$ whenever $\scr{P}_K$ contains an almost simple factor of type $A_n$, for some $n\geqslant 2$, any generically trivial $\scr{P}$-torsor over $\ca{O}_K$ is trivial, i.e., \bud \ker(H^1(\ca{O}_K,\scr{P})\to H^1(K,\scr{P}_K))=\{\ast\}.\eud
\eut
Our method, outlined below, proceeds by reducing Theorem~\ref{thm:A} to the previously discussed case when $\scr{P}$ is reductive. To carry out this reduction, we prove a fundamental structural result for all parahoric groups, namely, Theorem~\ref{thm:B}, whose statement will clarify our assumptions in Theorem~\ref{thm:A}. Roughly speaking, Theorem~\ref{thm:B} asserts that every parahoric group scheme $\scr{P}$ becomes reductive after a (possibly wildly ramified) extension of the base field $L/K$. More precisely, it shows that any parahoric group admits a \textit{reductification} after passing to a suitable extension $L$ of $K$. The notion of reductification, which is introduced below, is therefore a central concept in this article.
\bdf[Reductification]\label{defn:redutification} 
\bun
\item[$\circ$] A parahoric group scheme $\scr{P}$ over a henselian discrete valuation ring $\ca{O}_K$ with field of fractions $K$ is called \textit{reductifiable} if there exist a finite Galois extension $L/K$ with ring of integers $\ca{O}_L$ and Galois group $\Gamma$ and a $\Gamma$-equivariant reductive $\ca{O}_L$-model $\scr{G}$ of $\scr{P}\times_{\spec(K)}\spec(L)$, which has the canonical semilinear $\Gamma$-action induced by the action on the right factor, such that there is an isomorphism of $\ca{O}_K$-group schemes \bd\label{diag:intro:smooth_invariant_pushforward}\scr{P}\cong\fr{R}_{\ca{O}_L/\ca{O}_K}^{\mathit{\Gamma}}(\scr{G})_{\sm},\tag{$\ast$}\ed where $\fr{R}_{\ca{O}_L/\ca{O}_K}^{\mathit{\Gamma}}(\scr{G})$ denotes the closed subscheme parametrising the $\Gamma$-invariant sections (see \ccite{edixhoven_neron_model_tame_ramification}{\textsection3}) of the Weil restriction $\fr{R}_{\ca{O}_L/\ca{O}_K}(\scr{G})$ of $\scr{G}$ along the extension $\ca{O}_K\to\ca{O}_L$ and $(-)_{\sm}$ denotes the group smoothening functor (see Definition~\ref{def:dilatations}). 
We note that the $\Gamma$-action on $\fr{R}_{\ca{O}_L/\ca{O}_K}(\scr{G})$, representing $S\mapsto\scr{G}(S\times_{\spec(\ca{O}_K)}\spec(\ca{O}_L))$, is naturally induced by the $\Gamma$-action on $S\times_{\spec(\ca{O}_K)}\spec(\ca{O}_L)$ via the right factor for each $\ca{O}_K$-scheme $S$. If there is such an isomorphism \eqref{diag:intro:smooth_invariant_pushforward}, we will say that $(L,\scr{G})$ is a \textit{reductification} of $\scr{P}$.
\item[$\circ$] A parahoric group $\scr{G}$ over $\ca{O}_K$ is called \textit{tamely reductifiable} if it has a \textit{tame reductification}, i.e., a reductification $(L,\scr{G})$ where $L/K$ is tamely ramified.
In this case, post-composing with the group smoothening functor is unnecessary (see \ccite{edixhoven_neron_model_tame_ramification}{Proposition~3.4}) and we have an isomorphism of $\ca{O}_K$-group schemes \bud\scr{P}\cong\fr{R}_{\ca{O}_L/\ca{O}_K}^{\mathit{\Gamma}}(\scr{G}).\eud
\eun
\edf
The case where the generic fibre $\scr{P}_K$ is a torus is already exemplary and of independent interest. In this situation, the parahoric model $\scr{P}$ coincides with the identity component of the \textit{ft Néron model} of $\scr{P}_K$, uniquely characterised by the condition that $\scr{P}(\ca{O}_{\breve{K}})\subseteq\scr{P}_K(\breve{K})$ is the Iwahori subgroup (see \cite{neron_models}), where $\breve{K}$ is a maximal unramified extension of $K$. In particular, the torus $\scr{P}_K$ admits a reductification over its splitting field (see Remark~\ref{rem:torus_reductification}), which might be wildly ramified. This demonstrates the existence of reductifiable, but not tamely reductifiable, parahoric group schemes (see Example~\ref{exmp:edixhoven}). 
\par A second basic example occurs when $\scr{P}_K\cong\SL_{2,K}$, for which the only non-reductive parahoric model $\scr{P}$ is the Iwahori subgroup. We examine this case in Example~\ref{exmp:iwahori} and in Remark~\ref{rem:iwahori_not_tame}.
\par The main result of this article, stated below, generalises both of these situations discussed above.
\but{B}[Theorem~\ref{thm:reductification_parahoric} and Corollary~\ref{cor:tame_reductification_parahoric}]\label{thm:B}
	Every parahoric group scheme $\scr{P}$ over a henselian discrete valuation ring with a perfect residue field of characteristic $p$ is reductifiable. Moreover, this $\scr{P}$ is tamely reductifiable provided that its generic fibre $G$ is either a simply-connected, or adjoint, semisimple group scheme, and that $p$ satisfies the following conditions: $p\neq 2,3,5$ and $p\nmid (n+1)$ whenever $G$ contains an almost simple factor of type $A_n$, for some $n\geqslant 2$.
\eut
\br
Our proof strategy relies crucially on the techniques developed in \cite{kaletha_prasad}, where the residue field is assumed to be perfect by default (cf.~Remark~\ref{rem:perfect_residue_field1}). This assumption is harmless for our purposes, since the most interesting instance of Theorem~\ref{thm:B} is when the base ring is a finite extension of $\bb{Z}_p$ or $\bb{F}_p[[t]]$, in which case the residue field is even finite. Nevertheless, because the original source \cite{bruhat_tits_buildings_II} does impose this hypothesis, it should be possible to extend our result to a more general setting.
\er
The idea of reductification considered here originates in the work of \citeauthor{balaji_seshadri} \cite{balaji_seshadri}, who established the special case of Theorem~\ref{thm:B} where $p=0$ and $G$ is defined over $\bb{Z}$. Their work was generalised in \cite{damioli_thesis, damiolini_hong_local_types_equivariant_bundles}. In the tamely ramified setting, To prove Theorem~\ref{thm:B}, however, we follow the techniques of \citeauthor{pappas-rapoport_parahoric_grothendieck_serre} in \cite{pappas-rapoport_parahoric_grothendieck_serre}, who proved the second claim concerning tame reductification under the assumption that $G$ is simply-connected and almost simple (a related result in the tamely ramified setting was proved in \cite{chernousov_gille_pianzola_parahoric}).
\par Our main contribution is the treatment of the wildly ramified case, which becomes necessary whenever $G$ is a reductive group scheme whose central torus is not split by a tame extension, or when the residue characteristic is small enough. With this assumption, the situation is more delicate as counterexamples to tame reductification already occur; for instance, in the case of the Iwahori subgroup of $\SL_2$ over $\bb{Z}_2$ (see Remark~\ref{rem:iwahori_not_tame}). Consequently, post-composing with the group smoothening functor in the isomorphism \eqref{diag:intro:smooth_invariant_pushforward} is essential (see Example~\ref{exmp:edixhoven}). Group smoothening has previously been used to study parahoric models (see \cite{chai_yu_neron_model_tori, yu_smooth_model_Concave_functions}).
\par We expect Theorem~\ref{thm:B} to have applications beyond Question~\ref{quest:bayer_fluckigar_first}, for instance, in the study of moduli stacks of bundles over curves (see, for example, Balaji--Pandey \cite{balaji_pandey_conformal_blocks}, Damiolini \cite{damioli_thesis}, Damiolini--Hong \cite{damiolini_hong_local_types_equivariant_bundles}, and Heinloth \cite{heinloth_uniformisation}). In particular, it may be possible to extend the results in these papers to the wildly ramified setting.
\par We conclude this introduction with the promised sketch, followed by an outline of the article.
\subsubsection*{Overview of our proof of Theorem~\ref{thm:A}}
Adopting the notation of Theorem~\ref{thm:A}, Corollary~\ref{cor:tame_reductification_parahoric} shows that $\scr{P}$ admits a tame reductification $(L,\scr{G})$, i.e., there exists a reductive model $\scr{G}$ over a tamely ramified extension $L/K$ together with an isomorphism
\bud\scr{P} \cong \fr{R}_{\ca{O}_L/\ca{O}_K}^{\mathit{\Gamma}}(\scr{G}).\eud
By \ccite{pappas-rapoport_parahoric_grothendieck_serre}{Proposition~3.2} (cf.~\cite{damioli_thesis}), the $\scr{P}$-torsors (resp., with a generic trivialisation) over $\ca{O}_K$ correspond to $\scr{G}$-torsors over $\ca{O}_L$ with a compatible $\Gamma$-semilinear action (resp., with a $\Gamma$-equivariant generic trivialisation). These $(\Gamma,\scr{G})$-bundles, studied in \cite{balaji_seshadri} and \ccite{pappas-rapoport_parahoric_grothendieck_serre}{\textsection3.1} (cf.~\cite{damioli_thesis, damiolini_hong_local_types_equivariant_bundles}), are equivalent to torsors under the \textit{stacky reductive group} $[\scr{G}/\Gamma]$ over the \textit{stacky discrete valuation ring} $[\spec(\ca{O}_L)/\Gamma]$.
Our goal is a Levi reduction using the tameness of $L/K$. \par First, Lemma~\ref{lem:group_to_parabolic} produces a $\Gamma$-equivariant parabolic subgroup $\scr{Q} \subset \scr{G}$ over $\ca{O}_L$, so any generically trivial $[\scr{G}/\Gamma]$-torsor $\scr{E}$ admits a reduction to $[\scr{Q}/\Gamma]$. Second, Proposition~\ref{lem:parabolic_to_Levi} gives a $\Gamma$-equivariant Levi decomposition
$\scr{Q} \cong \scr{M} \rtimes \scr{U}$.
Using Serre vanishing for the tame stack $\mathbf{B}\Gamma$ (see \cite{abramovich_olsson_angelo_tame_stacks}), we show that $\scr{E}$ further reduces to a generically trivial $[\scr{M}/\Gamma]$-torsor.
\subsubsection*{Outline of this article}
In \textsection\ref{sec:prelims}, we recall the notion of parahoric group schemes following \cite{kaletha_prasad}; in \textsection\ref{sec:reductification}, we establish the first claim of Theorem~\ref{thm:B}; finally, in \textsection\ref{sec:gs_serre}, we prove the second claim of Theorem~\ref{thm:B}, thereby deducing Theorem~\ref{thm:A}.
\subsection*{Notations and Conventions}
By a \textit{discretely valued field}, we mean a field $K$ equipped with a discrete valuation. The ring of integers in such a field will be denoted $\ca{O}_K$, unless otherwise specified. This $K$ will be called \textit{henselian} if $\ca{O}_K$ is henselian. By \text{residue field} of $K$, we shall mean the residue field of $\ca{O}_K$, which is assumed to be perfect. The maximal unramified extension of $K$ will be denoted by $\breve{K}$.
\section{Kottwitz homomorphism and the Bruhat--Tits building.}\label{sec:prelims}
In this section, following \cite{kaletha_prasad}, our goal is to define the notion of parahoric group schemes in Definition~\ref{defn:parahoric}. As highlighted in the introduction, these are certain integral models of reductive group schemes $G$ defined over henselian discretely valued fields, characterised by the property that their group of integer points coincide with the connected stabilisers of points in the Bruhat--Tits buildings associated to $G$. We review necessary properties of the Bruhat--Tits buildings in \textsection\ref{mo:bruhat-tits_buildings}. Moreover, since Theorem~\ref{thm:A} concerns non-simply-connected groups schemes, the notion of the Kottwitz homomorphism recalled \textsection\ref{mo:kottwitz_homomorphism} is essential in Definition~\ref{defn:parahoric}.
\par We begin with the following property of fields, our adaptation of \ccite{serre_galois_cohomology}{Chapter~2, \textsection3.1} and \ccite{kaletha_prasad}{\textsection2.3}; this is satisfied by the maximal unramified extension of all the base fields in consideration in this article.
\bdf[Cohomological dimension]\label{defn:cd_1}
	A field $F$ is said to have \textit{cohomological dimension $\leqslant 1$}, denoted $\cd(F)\leqslant 1$, if for any algebraic extension $L/F$ the following equivalent conditions hold.
	\bn[(1)]
		\item For any finite Galois extension $L'/L$, the associated norm morphism $N_{L'/L}\colon \bb{G}_m(L')\to\bb{G}_m(L)$ is surjective.
		\item For any $L$-torus $T$, we have $H^1(L,T)=0$.
	\en
\edf
The equivalence of the two conditions in Definition~\ref{defn:cd_1} follows from Steinberg's theorem (see \ccite{kaletha_prasad}{Theorem~2.3.3}).
\br\label{rem:unramified_extension_cd_1}
	\bn[(1)]
		\item\label{rem:unramified_extension_cd_1:pt:1} More generally, Steinberg's theorem proves that for a field $F$ to have $\cd(F)\leqslant 1$ it is necessary and sufficient that for any algebraic extension $L/F$, every reductive $L$-group $G$ satisfies $H^1(L,G)=\{\ast\}$.
		\item\label{rem:unramified_extension_cd_1:pt:2} Furthermore, the same theorem guarantees that over a field $F$ satisfying $\cd(F)\leqslant 1$, any reductive group $G$ admits a Borel $F$-subgroup, in other words, any such $G$ is $F$-\textit{quasi-split}.
		\item\label{rem:unramified_extension_cd_1:pt:3} The construction of the Bruhat--Tits building associated to a reductive group $G$ defined over a henselian discretely valued field $K$ in \cite{bruhat_tits_buildings_II} requires that $G$ be quasi-split over the maximal unramified extension $\breve{K}$ of $K$. By \eqref{rem:unramified_extension_cd_1:pt:2}, this is satisfied as long as $\cd(\breve{K})\leqslant 1$, which follows from our assumption throughout this article that $K$ is henselian with a perfect residue field (see \ccite{kaletha_prasad}{Corollary~2.3.7}).
	\en
\er
Let us introduce the following notations used in rest of this section. Given a torus $T$ over a field $F$ with a separable closure $F^{\sep}$, the \textit{character lattice} (resp., \textit{cocharacter lattice}) is denoted by $\bb{X}^{\ast}(T)\colonequals\Hom_{F^{\sep}}(T,\bb{G}_{m})$ (resp., $\bb{X}_{\ast}(T)\colonequals\Hom_{F^{\sep}}(\bb{G}_{m}, T)$). Either of these objects is naturally equipped with an action of the absolute Galois group $\Gal(F^{\sep}/F)$.
\par The definition in \textsection\ref{mo:kottwitz_homomorphism} of the Kottwitz homomorphism associated to a reductive group is rather involved, relying on various notions that we need to introduce before. We start with the first one, called a \textit{$z$-extension} (synonymous to a \textit{$t$-resolution} in \cite{borovois_gonzalez_algebraic_fundamental_group} and which includes the notion of a \textit{coflasque resolution} in \cite{ct_borovoi_fundamental_group}), introduced by Langlands and Ono and effectively utilised by Kottwitz. To do so, we require first to introduce the following concept, which is called differently throughout the literature, e.g., \textit{quasi-trivial} in \cite{ct-sansuc_flasque_tori}. However, we adopt the name chosen by \citeauthor{kaletha_prasad} in \ccite{kaletha_prasad}{Definition~2.5.1}.
\bdf[Induced torus]\label{defn:induced_torus} 
	A torus $T$ over a field $F$ is called \textit{induced} if there exists a Galois splitting field extension $F'/F$ such that the lattice $\bb{X}^{\ast}(T)$, equivalently, the lattice $\bb{X}_{\ast}(T)$, has a $\Gal(F'/F)$-stable $\bb{Z}$-basis.
\edf
\br\label{rem:induced_torus}
	Let $T$ be an induced torus over a field $F$.
	\bn[(i)]
		\item\label{rem:induced_torus:pt:1} It follows from Definition~\ref{defn:induced_torus} that $T$ can be written as a finite product of $\Res_{L/F}(\bb{G}_m)$ where each $L$ is a certain subfield of a chosen Galois splitting field extension $F'/F$, i.e., an extension $F'/F$ for which $\bb{X}^{\ast}(T)$, equivalently, $\bb{X}_{\ast}(T)$, has a $\Gal(F'/F)$-equivariant $\bb{Z}$-basis. Indeed, this follows from the equivalence between the category of $F'$-split $F$-tori and that of finite free $\bb{Z}$-modules endowed with a $\Gal(F'/F)$-action \ccite{sga3ii}{exposé~X, proposition~1.1}.
		\item In particular, \eqref{rem:induced_torus:pt:1} shows that given any Galois splitting field extension $F''/F$, the lattice $\bb{X}^{\ast}(T)$, equivalently, the lattice $\bb{X}_{\ast}(T)$, has a $\Gal(F''/F)$-stable $\bb{Z}$-basis.
		\item The maximal torus of quasi-split semisimple group scheme that is either adjoint or simply-connected is induced (see \ccite{bruhat_tits_buildings_II}{proposition~4.4.16}).
	\en
\er
\bdf[$Z$-extension]\label{defn:z-extension} Given a reductive group scheme $G$ over a field $F$, a \textit{$z$-extension} (resp., a \textit{$\mu$-extension}) is an exact sequence \bud 1\to Z\to\widetilde{G}\to G\to 1,\eud where $Z$ is an induced $F$-torus (resp., where $Z$ is a finite $F$-group scheme of multiplicative type) and $\tilde{G}$ is a reductive $F$-group scheme whose derived subgroup is simply-connected.
\edf
\br\label{rem:z-extension}
	\bn[(a)] 
		\item\label{rem:z-extension:pt:1} Given either a $z$-extension or a $\mu$-extension $1\to Z\to\tilde{G}\to G\to 1$, since $\widetilde{G}$ is connected, the normal subgroup $Z\subset\widetilde{G}$ of multiplicative type is necessarily central. This can be proved by \ccite{milne_algebraic_groups}{Corollary~12.38} or by noting that the automorphism group of a split group of multiplicative type is discrete (see \ccite{sga3ii}{exposé~X, proposition~1.1}); for example, in the case of a split torus of rank $r$ this group is $\GL_r(\bb{Z})$.
		\item\label{rem:z-extension:pt:2} Given a reductive $F$-group scheme $G$, we recall the construction of the $z$-extension from \ccite{kaletha_prasad}{Construction~11.4.1}. First, we know that the identity component $T$ of the centre of $G$ is a torus (\ccite{sga3ii}{exposé~XII, proposition~4.11}). Second, we let $G_{\SC}$ be the simply-connected cover of the derived subgroup of $G$ (see \ccite{conrad_reductive_group_schemes}{Exercise~1.6.13}) and $\mu$ be the finite multiplicative group that is the kernel of the surjective isogeny (see \ccite{borel_algebraic_groups}{\textsection14.2}) \bd\label{rem:mu_extension} T\times G_{\SC}\twoheadrightarrow G.\ed This furnishes a $\mu$-extension $1\to\mu\to T\times G_{\SC}\to G\to 1$. 
		\par To produce a $z$-extension, we now suppose that $\mu$ embeds inside a torus $Z$, the existence of which is can be proved as in \ccite{kaletha_prasad}{Construction~11.4.1}. Thereafter, we obtain a $z$-extension of $G$ by taking \bud\tilde{G}\colonequals (Z\times T\times G_{\SC})/\mu.\eud
	\en
\er
The following, recalled from \ccite{kaletha_prasad}{\textsection2.5(b)}, is an additional ingredient in the definition of the Kottwitz homomorphism.
\bmo[Valuation Homomorphism]\label{mo:valuation_homomorphism} Given a torus $T$ over a henselian discretely valued field $(F,|\cdot|_F)$ with a perfect residue field and an absolute Galois group $\fr{G}$, we have a functorial homomorphism $\omega_{T,F}\colon T(F)\to\Hom_{\bb{Z}}(\bb{X}^{\ast}(T)^{\fr{G}}, \bb{Z})$ given by $t\mapsto (\chi\mapsto -|\chi(t)|_{F})$. When $T$ is $F$-split, using the isomorphism $\Hom_{\bb{Z}}(\bb{X}^{\ast}(T),\bb{Z})\cong \bb{X}_{\ast}(T)$, we obtain the following surjective morphism (\ccite{kaletha_prasad}{Lemma~2.5.7}) called the \textit{valuation homomorphism} \bud v_{T,F}\colon T(F)\to \bb{X}_{\ast}(T).\eud
\emo
Before introducing the Kottwitz homomorphism in \textsection\ref{mo:kottwitz_homomorphism}, the final necessary concept (known also in the literature as the \textit{algebraic fundamental group}) is the following, recalled from \ccite{kaletha_prasad}{\textsection11.3}.
\bmo[Borovoi fundamental group]\label{mo:borovoi_fundamental_group}
Let $F$ be a field with a separable closure $F^{\sep}$ and $G$ be a reductive $F$-group scheme. Suppose that $T\subseteq G$ is a maximal $F$-torus and $G_{\SC}$ is the simply-connected cover of the derived subgroup of $G$. Let $T_{\SC}\subseteq G_{\SC}$ be a maximal $F$-torus that lies in the preimage of $T$. The cocharacter lattice of $T_{\SC}$ can then be identified with the sublattice $\bb{X}_{\ast}(T_{\SC})\subseteq \bb{X}_{\ast}(T)$ generated by the coroots of the split reductive $F^{\sep}$-group $(G_{F^{\sep}},T_{F^{\sep}})$ (see \ccite{ces_torsors}{Proposition~A.3.4}). Then, the \textit{Borovoi fundamental group} of $G$ is defined to be the following finitely generated abelian group with a natural $\Gal(F^{\sep}/F)$-action \bud\pi^
B_1(G)\colonequals \bb{X}_{\ast}(T)/\bb{X}_{\ast}(T_{\SC}).\eud 
This construction is independent of the choice of the maximal tori $T\subset G$ (see \ccite{borovois_gonzalez_algebraic_fundamental_group}{\textsection2} or \ccite{ct_borovoi_fundamental_group}{proposition~A.2}); and moreover, the assignment $G\mapsto\pi^B_1(G)$ depends functorially in $G$ (see \ccite{kaletha_prasad}{Lemma~11.3.6}).
\par We note that, by construction, if $G=T$ is a torus, there is an equality $\pi^B_1(T)=\bb{X}_{\ast}(T)$. More generally, let $G_{\der}\subseteq G$ be the derived subgroup and $G_{\text{ab}}\colonequals G/G_{\der}$ be the maximal abelian quotient of $G$. In the case where $G_{\der}$ is simply-connected, using the isomorphism between the character and the cocharacter lattices, we have a canonical equality \bd\label{iso:derived_simply_connected_fundamental_group}\pi_1^B(G)=\pi^B_1(G_{\text{ab}}).\ed
\emo
We will use the following notation in the rest of the section.
\bmo[Setup]\label{mo:setup}~
\bun
\item[$\diamond$] Let $K$ be a henselian discretely valued field with a perfect residue field. 
\item[$\diamond$] Let $\breve{K}$ be a maximal unramified extension of $K$ contained in a separable closure $K^{\sep}$ of $K$. 
\item[$\diamond$] Let $I\colonequals\Gal(K^{\sep}/\breve{K})$ be the inertia subgroup of the absolute Galois group $\Gamma\colonequals\Gal(K^{\sep}/K)$ with cokernel $\Xi\colonequals\Gamma/I\cong\Gal(\breve{K}/K)$.
\eun  \emo
\bmo[Kottwitz homomorphism]\label{mo:kottwitz_homomorphism} Using notations introduced in \textsection\ref{mo:setup}, given a reductive $K$-group scheme, we recall the definition of the \textit{Kottwitz homomorphism}
\bud \kappa_G\colon G(\breve{K})\to\pi_1^B(G)_I, \eud from \ccite{kaletha_prasad}{\textsection11.5}. The construction, written below, shows that this morphism is surjective and $\Xi$-equivariant, simultaneously, which functorially agrees with a base field extension as well as with a morphism of reductive groups schemes. 

\bn[(i)]
	\item\label{mo:kottwitz_homomorphism:pt:1} Let $T$ be a $K$-torus, and, for each finite Galois extension $L/\breve{K}$, let $N_{L/\breve{K}}\colon T(L)\to T(\breve{K})$ be the \textit{norm morphism} associated to $\breve{K}\to L$ (see, for example, \ccite{phd-thesis}{\textsection5.5}). We note that this morphism is surjective due to our assumption that $\cd(\breve{K})\leqslant 1$ (see Remark~\ref{rem:unramified_extension_cd_1}\eqref{rem:unramified_extension_cd_1:pt:1}). The Kottwitz homomorphism for $T$ is therefore defined to be the unique morphism that makes the following diagram commute for any finite Galois extension $L/\breve{K}$ splitting $T$ \bud T(L)\arrow{d}[swap]{N_{L/\breve{K}}}\arrow{r}{v_{T,L}}[swap]{\eqref{mo:valuation_homomorphism}}& \bb{X}_{\ast}(T)\arrow[d]\\ T(\breve{K})\arrow[r,"\kappa_T"]&\bb{X}_{\ast}(T)_I,\eud where we have made the identification $\bb{X}_{\ast}(T)=\pi^B_1(T)$ and where the right vertical arrow is the canonical quotient morphism. The $\Xi$-equivariance of $\kappa_T$ follows automatically from the uniqueness and the $\Xi$-equivariance of $v_{T,L}$.
	\item\label{mo:kottwitz_homomorphism:pt:2} In the case where the derived subgroup $G_{\der}$ of $G$ is simply-connected, the Kottwitz homomorphism is defined so that the diagram below commutes \bud G(\breve{K})\arrow[r,"\kappa_G"]\arrow[d]&\pi_1^B(G)_I\arrow[d,equal]\\G_{\mathrm{ab}}(\breve{K})\arrow{r}{\eqref{mo:kottwitz_homomorphism:pt:1}}[swap]{\kappa_{G_{\mathrm{ab}}}}&\pi^B_1(G_{\mathrm{ab}})_I, \eud where $G_{\mathrm{ab}}\colonequals G/G_{\der}$ is a torus (see \ccite{conrad_reductive_group_schemes}{Theorem~5.3.1}) and the left equality is \eqref{iso:derived_simply_connected_fundamental_group}. The $\Xi$-equivariance is then a consequence of the $\Xi$-equivariance of $\kappa_{G_{\mathrm{ab}}}$ and the same of the left vertical arrow.
	\item In general, we choose a $z$-extension (for instance, the $z$-extension constructed in Remark~\ref{rem:z-extension}\eqref{rem:z-extension:pt:2}) \bud 1\to Z\to\widetilde{G}\to G\to 1\eud and define the Kottwitz homomorphism such that the following digram commutes \bud \widetilde{G}(\breve{K})\arrow[d]\arrow{r}{\kappa_{\widetilde{G}}}[swap]{\eqref{mo:kottwitz_homomorphism:pt:2}}&\pi_1^B(\widetilde{G})_I\arrow[d]\\ G(\breve{K})\arrow{r}[swap]{\kappa_G}&\pi_1^B(G)_I,\eud where left vertical arrow is surjective (since $\cd(\breve{K})\leqslant 1$ by Remark~\ref{rem:unramified_extension_cd_1}\eqref{rem:unramified_extension_cd_1:pt:3}) and so is the right vertical morphism (see \ccite{borovois_gonzalez_algebraic_fundamental_group}{\textsection2}, cf.~\ccite{ct_borovoi_fundamental_group}{proposition~6.2}). The $\Xi$-equivariance of $\kappa_G$ is a direct consequence of the analogous property for $\kappa_{\widetilde{G}}$ and the right vertical arrow. Moreover, the independence of the construction from the chosen $z$-extension (see \ccite{kaletha_prasad}{Lemma~11.5.1}) guarantees that the resulting definition is well-defined.
\en
Taking $\Xi$-invariants then reduces $\kappa_G$ to the following homomorphism, also called the Kottwitz homomorphism by abuse of terminology, \bud \kappa_G|_{G(K)}\colon G(K)\to(\pi_1^B(G)_I)^{\Xi}. \eud
As an outcome of the above discussion, we define the following distinguished subgroups that are essential to characterise parahoric group schemes associated to non-simply-connected groups. We let \bd\label{defn:G(K)^0} G(\breve{K})^0\colonequals\ker(\kappa_G)~~\text{ and }~~G(K)^0\colonequals\ker(\kappa_G|_{G(K)}).\ed 
Let us note that \ccite{kaletha_prasad}{Proposition~11.5.4} shows that there is an equality $G(K)^0=G(\breve{K})^0\cap G(K)$.
\emo
\br[Kottwitz homomorphism under base extension]\label{rem:functoriality_of_Kottwitz}
	By construction, for any finite extension of discretely valued fields $L/K$, there is a commutative diagram \bud G(\breve{K})\arrow[d]\arrow{r}{\kappa_{G}}[swap]{}&\pi_1^B(G)_{I}\\ G(\breve{L})\arrow{r}[swap]{\kappa_{G_L}}&\pi_1^B(G_L)_{I_L}\arrow[u],\eud where $\breve{L}$ is a maximal unramified extension of $L$ and $I_L$ is the inertia subgroup of $L$. In particular, this shows that we have an inclusion \bud G(\breve{L})^0\cap G(\breve{K})\subseteq G(\breve{K})^0.\eud
\er
\bmo[Bruhat--Tits building]\label{mo:bruhat-tits_buildings} Adopting notations in \textsection\ref{mo:setup}, given a reductive $K$-group scheme $G$, there exist two (restricted) \textit{abstract buildings} in the sense of \ccite{kaletha_prasad}{Definition~1.5.5}, i.e., two polysimplicial complexes satisfying a list of properties enunciated in loc.~cit. Namely,
\bun
	\item the \textit{Bruhat--Tits building} $\scr{B}(G,\breve{K})$ of $G$ over $\breve{K}$ defined in \ccite{kaletha_prasad}{Definition~7.6.1}, and
	\item the \textit{Bruhat--Tits building} $\scr{B}(G,K)\colonequals\scr{B}(G,\breve{K})^{\Xi}$ of $G$ over $K$ defined in \ccite{kaletha_prasad}{Definition~9.2.8}.
\eun
The building $\scr{B}(G,\breve{K})$ comes equipped with an action of $G(\breve{K})\rtimes\Gamma$; consequently, $\scr{B}(G,K)$ comes equipped with an action of $G(K)$. Additionally, the construction respects finite extensions of base fields. Indeed, given a finite extension of discretely valued fields $L/K$, by \ccite{kaletha_prasad}{\textsection7.9.2}, there is a natural injection, equivariant with respect to $G(\breve{K})\subseteq G(\breve{L})$, \begin{equation}\label{diag:building_extension_base_field}
	\begin{tikzcd}\scr{B}(G,\breve{K})\arrow[r,hook]&\scr{B}(G,\breve{L}).\end{tikzcd}
\end{equation}
If the extension $L/K$ is Galois (resp., Galois and tamely ramified), the image of the above displayed injection lies inside (resp., equals) the $\Gal(L/K)$-invariant subset of the target (see \ccite{kaletha_prasad}{Remark~7.9.3}).
Furthermore, by \ccite{kaletha_prasad}{Proposition~14.1.1}, given a surjective morphism $G\twoheadrightarrow H$ of reductive $K$-groups, there exists a surjective morphism \bd\label{diag:building_quotient_map}\scr{B}(G,\breve{K})\twoheadrightarrow\scr{B}(H,\breve{K})\ed that is equivariant with respect to $G(\breve{K})\twoheadrightarrow H(\breve{K})$.
\emo
\br\label{rem:perfect_residue_field1}
	The perfectness of the residue field of $K$ in \textsection\ref{mo:setup} is assumed solely to ensure that $G_{\breve{K}}$ is quasi-split (see Remark~\ref{rem:unramified_extension_cd_1}\eqref{rem:unramified_extension_cd_1:pt:2}), that is a requirement in the construction in \ccite{kaletha_prasad}{Chapter~9}. This hypothesis therefore becomes unnecessary if one directly assumes that $G$ is quasi-split over $\breve{K}$.
\er
Finally, we turn to the definition of parahoric group schemes.
\bdf[Parahoric group]\label{defn:parahoric}
Let us use the notations from \textsection\ref{mo:setup}.
\bn[(a)] 
\item\label{defn:parahoric_subgroup} Given a reductive $\breve{K}$-group $G$, a subgroup $\ca{P}\subseteq G(\breve{K})$ is said to be \textit{parahoric} if there exists a point $x\in\scr{B}(G,\breve{K})$ such that $\ca{P}=G(\breve{K})^0\cap\Stab_{G(\breve{K})}(x)$ (see \eqref{defn:G(K)^0}).
\item\label{defn:parahoric:1} An  affine, smooth, fibrewise connected $\ca{O}$-group scheme $\scr{P}$ is called \textit{parahoric} if its generic fibre $\scr{P}_K$ is a reductive $K$-group scheme and $\scr{P}(\ca{O}_{\breve{K}})\subseteq\scr{P}_K(\breve{K})$ is a parahoric subgroup in the sense of \eqref{defn:parahoric_subgroup}. In this holds, this $\scr{P}$ is then called a \textit{parahoric model} of $\scr{P}_K$.
\item\label{defn:parahoric:2} More generally, an  affine, smooth, fibrewise connected group scheme $\ca{G}$ over a Dedekind ring $R$ with perfect residue fields is called \textit{parahoric} if for each maximal ideal $\fr{m}\subset R$ the base change of $\ca{G}$ to the henselisation of $R$ at $\fr{m}$ is a parahoric group in the sense of \eqref{defn:parahoric:1}.
\en
\edf
\section{Parahoric groups as smooth invariant pushforwards of reductive groups.}\label{sec:reductification}
In this section, our aim is to prove Theorem~\ref{thm:reductification_parahoric}, which establishes that every parahoric group scheme admits a (possibly wild) reductification, thereby establishing the first claim in Theorem~\ref{thm:B}. To do so, we follow the techniques of Pappas and Rapoport \cite{pappas-rapoport_parahoric_grothendieck_serre} (which are similar in spirit to those of Balaji and Seshadri in \cite{balaji_seshadri}). An important step in the proof of this theorem is Proposition~\ref{prop:larsen_hyperspeciality}, which is essentially \ccite{pappas-rapoport_parahoric_grothendieck_serre}{Proposition~2.7}. We follow the techniques in loc.~cit. and deduce it as a consequence of Proposition~\ref{lem:larsen_hyperspecial}, a result that has appeared in the work of Larsen \cite{larsen_maximality} and Gille \cite{gille_torseurs_sur_la_droite_affine} (cf.~Cotner \cite{cotner_parahoric}).
\par As mentioned in the introduction, the primary novelty of this article is treating the case in which the parahoric group scheme has a wild reductification (see Theorem~\ref{thm:reductification_parahoric}). In this situation, the invariant pushforwards of group schemes are not automatically smooth, as noted in Example~\ref{exmp:edixhoven}. In fact, we need to enforce group smoothening separately, and the resulting object is then called the \textit{smooth invariant pushforward}. The idea of using group smoothening to construct parahoric models has previously appeared in \cite{chai_yu_neron_model_tori, yu_smooth_model_Concave_functions}. The notion of group smoothening is recalled from \ccite{neron_models}{\textsection7.1} in Definition~\ref{defn:smoothening}.
\par We begin with the following concept from \ccite{neron_models}{\textsection3.2}, that is required to study properties of group smoothening in Remark~\ref{rem:smoothening}.
\bmo[Dilatation]\label{mo:dilatation}
\emo
Given a scheme $X$ over a discrete valuation ring $\ca{O}$, let \bud X_{\mathrm{fl}}\hookrightarrow X\eud be the largest flat $\ca{O}$-subscheme, which is given by the vanishing of the quasi-coherent sheaf of ideals generated by the $\pi^{\infty}$-torsion elements. This subscheme satisfies the following universal property: any $\ca{O}$-morphism $Y\to X$ from a $\ca{O}$-flat (equivalently, $\ca{O}$-torsion free) $\ca{O}$-scheme factors uniquely through through $X_{\mathrm{fl}}\hookrightarrow X$. 
\bdf[Dilatation]\label{def:dilatations}
	Given a scheme $X$ over a discrete valuation ring $\ca{O}$, the \textit{dilatation} (synonymously, the \textit{affine blowup}) of $X$ along a closed subscheme $Z$ of the $\ca{O}$-special fibre of $X$ is the terminal object \bud\Bl^{\aff}_Z(X)\to X\eud in the category of $\ca{O}$-morphisms $Y\to X$ whose restriction to the $\ca{O}$-special fibre factors through $Z$ and whose source $Y$ is a flat $\ca{O}$-scheme.
\edf
We summarise \ccite{neron_models}{\textsection3.2, Propositions~1-3} in the following remark.
\br\label{rem:dilatations}
	Suppose $\ca{O}$ is a discrete valuation ring with a uniformiser $\pi\in\ca{O}$. Let $X$ be an $\ca{O}$-scheme and $Z$ be a closed subscheme of the $\ca{O}$-special fibre $X/(\pi)$ of $X$.
	\bn[(a)]
		\item\label{pt:a:dilatations} The dilatation $\Bl^{\aff}_Z(X)\to X$ exists and factors uniquely through the largest flat $\ca{O}$-subscheme $X_{\mathrm{fl}}\hookrightarrow X$. To construct this dilatation, by uniqueness, it suffices to assume that $X$, and hence, $Z$ is affine. Let $X=\spec(A)$ and $Z=V(I)$ for an ideal $I\subset A$. Since $A$ is noetherian, we have $I=(\pi,g_1,\ldots,g_n)$ for some elements $g_i\in A$. Then $\Bl^{\aff}_Z(X)\subsetneq\Bl_Z(X)$ corresponds to the open locus where this $I$ is generated by $\pi$. We verify then that $\Bl^{\aff}_{Z}(X)=(\spec(A'))_{\mathrm{fl}}$ where $A'=A[g_1/\pi,\ldots,g_n/\pi]$.
		\item Given an $\ca{O}$-morphism $f\colon X\to X'$ and a closed subscheme $Z'\subseteq X'/(\pi)$ containing $f(Z)$, there is a canonical morphism $\Bl^{\aff}_Z(X)\to\Bl^{\aff}_{Z'}(X')$.
	\en
\er
\bdf[Smoothening]\label{defn:smoothening}
	A \textit{smoothening} of a generically smooth, finite type group scheme $\scr{G}$ over a discrete valuation ring $\ca{O}$ is the terminal object \bud\scr{G}_{\sm}\to\scr{G}\eud in the category of $\ca{O}$-morphisms $\scr{X}\to\scr{G}$ whose source is a smooth $\ca{O}$-scheme $\scr{X}$ of finite type.
\edf
\br\label{rem:smoothening}
	Let $\scr{G}$ be a generically smooth, group scheme over a discrete valuation ring $\ca{O}$.
	\bn[(a)]
	\item It follows from the definition that a smoothening $\scr{G}_{\sm}\to\scr{G}$ is unique if it exists. Furthermore, the smoothening induces a isomorphism \bd\label{iso:generic_Smoothening}\scr{G}_{\sm,K}\iso\scr{G}_K\ed of fibres over the fraction field $K$ of $\ca{O}$.
	\item Granted the existence of $\scr{G}_{\sm}$, the $\ca{O}$-scheme $\scr{G}_{\sm}\times\scr{G}_{\sm}$ is then also $\ca{O}$-smooth of finite type. Consequently, the uniqueness guarantees that this $\scr{G}_{\sm}$ obtains the structure of an $\ca{O}$-group scheme that is respected by the canonical morphism $\scr{G}_{\sm}\to\scr{G}$. In this case, the universality further guarantees that the generic isomorphism \eqref{iso:generic_Smoothening} is compatible with their group structures.
	\item It follows from Definition~\ref{defn:smoothening} that smoothening is functorial provided it existence, i.e., given a morphism of generically smooth, finite type $\ca{O}$-group schemes $\scr{G}'\to\scr{G}$ of finite type, there exists a commutative diagram of $\ca{O}$-group schemes \bud \scr{G}'_{\sm}\tir\arrow[d]&\scr{G}'\arrow[d]\\\scr{G}_{\sm}\tir&\scr{G}.\eud
	\item This $\scr{G}_{\sm}$ may be constructed by a finite sequence of dilatations of $\scr{G}$ along closed subgroup schemes of positive codimension of the $\ca{O}$-special fibre of $\scr{G}$. Indeed, by \ccite{neron_models}{\textsection7.1, Lemma~4}, each such dilatation reduces \textit{Néron's measure for the defect of smoothness}, which is a non-negative integer defined in \ccite{neron_models}{\textsection3.3}.
	\item In particular, the above and Remark~\ref{rem:dilatations}\eqref{pt:a:dilatations} ensures that if $\scr{G}$ is an affine $\ca{O}$-scheme, then so is $\scr{G}^{\sm}$. 
	\en
\er
The following is essentially contained in \ccite{neron_models}{\textsection7.1, Theorem~5}.
\bp\label{prop:smoothening_exists}
	Given a generically smooth, finite type group scheme $\scr{G}$ over a discrete valuation ring $\ca{O}$, the category of $\ca{O}$-group homomorphisms $\scr{H}\to\scr{G}$ whose source is a smooth $\ca{O}$-group scheme $\scr{H}$ of finite type has a terminal object $\scr{G}_{\sm}\to\scr{G}$. Moreover, over the strict henselisation $\ca{O}^{\sh}$ of $\ca{O}$, this terminal object satisfies the following bijections
	\bud\ast_0\colon\scr{G}_{\sm}(\ca{O}^{\sh})\iso\scr{G}(\ca{O}^{\sh})~~~~~~\text{ and }~~~~~~\ast_{1}\colon H^1_{\etale}(\ca{O},\scr{G}_{\sm})\iso H^1_{\etale}(\ca{O},\scr{G}).\eud 
\ep
In the rest of the article, given a group scheme $G$ over a scheme $S$, the groupoid of $G$-torsors over $S$ that are trivialised étale-locally over $S$ will be denoted by $\mathbf{B}G(S)$.
\bs
	First, assuming that $\ast_0$ induces a bijection, we show that there is an equivalence \bd\label{diag:_smoothening_torsors_equivalence}\mathbf{B}\scr{G}_{\sm}(\ca{O})\iso\mathbf{B}\scr{G}(\ca{O}),\ed which ensures, by taking isomorphism classes of objects, that $\ast_1$ is a bijection. Thanks to the bijection $\ast_0$, in order to show \eqref{diag:_smoothening_torsors_equivalence}, it suffices to note that the $\scr{G}_{\sm}\to\scr{G}$ induces an isomorphism of sheaves on the small étale site over $\ca{O}$. 
\es
\br
	When $\scr{G}$ is not $\ca{O}$-smooth, diagram \eqref{diag:_smoothening_torsors_equivalence} precisely establishes an equivalence between torsors over $\ca{O}$ that are trivialised étale-locally over $\ca{O}$.
\er
To each simple reductive group over a $p$-adic field, there is an associated set of bad primes, denoted $S(G)$ by Serre in \ccite{serre_survery_galois_cohomology}{\textsection2.2}. Inspired by this and the work of Pappas and Rapoport in \cite{pappas-rapoport_parahoric_grothendieck_serre}, we introduce the following notion, which is required to state Proposition~\ref{prop:larsen_hyperspeciality} below. 
\bdf[Bad simple factor]\label{defn:Serre_number}
	Let $K$ be a henselian discretely valued field with a perfect residue field of characteristic $p$ and $H$ be a semisimple $K$-group scheme having a decomposition into absolutely simple $\breve{K}$-group schemes $H\times_{\spec(K)}{\spec(\breve{K})}=\prod_i H_i$ as in \ccite{sga3iii}{exposé~XXVI, corollaire~6.12}. For each $i$, this $H_i$ is called a \textit{bad simple factor} of $H$ if $$\begin{cases}p\mid 2(n+1)~~~\text{ and }H_i\text{ is of type }A_n,\text{ for some }n\geqslant 2;\\ p= 2~~~~~~~~~~~~~~\text{ and }H_i\text{ is of type }B_n, C_n, D_n\text{ (for }n\neq 4\text{), or }G_2;\\p= 2,3~~~~~~~~~~\text{ and }H_i\text{ is of type }D_4,F_4,E_6\text{ or }E_7\\ p= 2,3,5~~~~~~\text{ and }H_i\text{ if of type }E_8;\end{cases}$$
	If such a factor does not exist, then $H$ is said to have \textit{no bad simple factor}.
\edf
The following, which is the key ingredient in Proposition~\ref{prop:larsen_hyperspeciality}, has appeared in \ccite{larsen_maximality}{Lemma~2.4}, \ccite{gille_torseurs_sur_la_droite_affine}{\textsection2.2}, and \ccite{cotner_parahoric}{Theorem~4.3}).
	\bp\label{lem:larsen_hyperspecial}
		Let $K$ be a henselian discretely valued field with a perfect residue field and $H$ be a reductive $K$-group scheme. Given a barycentre of a facet $x\in\scr{B}(H,\breve{K})$, there exists a finite Galois extension $L/K$ of discretely valued fields splitting $H$ such that $x$ maps via \eqref{diag:building_extension_base_field} to a hyperspecial point in $\scr{B}(H,\breve{L})$. Furthermore, if $H$ is semisimple with no bad simple factor in the sense of Definition~\ref{defn:Serre_number}, this extension $L/K$ can be chosen to be tamely ramified.
	\ep
	\bs
		This essentially follows from the arguments of \ccite{larsen_maximality}{Lemma~2.4} (cf.~\ccite{cotner_parahoric}{Theorem~4.3}). We note that even though the statement of \ccite{larsen_maximality}{Lemma~2.4} requires $H$ to split over an unramified extension, this hypothesis can be removed. We sketch a modification of the argument in loc. cit. that works without this hypothesis.
		\par Let $F/K$ be a finite Galois extension splitting $H$ and let $x$ be the corresponding image in $\scr{B}(H,\breve{F})$. In the case where $H$ is semisimple and does not have any bad simple factor, we can choose $F/K$ to be tamely ramified. Indeed, in this case, the hypothesis ensures that the order of the finite group of automorphisms $\Aut(\scr{R})$ of the absolute affine Dynkin diagram (equivalently, the absolute extended Dynkin diagram) $\scr{R}$ is coprime to $p$. Therefore, since the base $H_{\breve{K}}$ is quasi-split (see Remark~\ref{rem:unramified_extension_cd_1}\eqref{rem:unramified_extension_cd_1:pt:3}), \bud\text{the associated absolute Galois representation $\Gal(K^{\sep}/\breve{K})\to\Aut(\scr{R})$ is tamely ramified.}\eud By \ccite{sga3iii}{exposé~XXVI, corollaire~6.12}, this $H_{\breve{K}}$ decomposes as a product of absolutely simple $\breve{K}$-group schemes, reducing us to the case where $H_{\breve{K}}$ is an absolutely simple $\breve{K}$-group scheme. Consequently, by the classification of absolutely simple $\breve{K}$-group scheme in \ccite{kaletha_prasad}{\textsection10.7}, it follows that $H_{\breve{K}}$ splits over a tamely ramified extension of $\breve{K}$. Fianlly, a limit argument shows that $H$ itself splits over a tamely ramified extension $F/K$.
		\par  Subsequently, as in the footnote in \ccite{cotner_parahoric}{Theorem~4.3}, possibly by replacing $F$ with a further finite Galois extension, we may assume that $x\in\scr{B}(H,\breve{F})$ is a vertex. Since $H_{\breve{F}}$ is split, by \ccite{kaletha_prasad}{Definition~7.11.1 and Remark~7.11.2}, for any extension $L/F$, the notion of special and hyperspecial points in $\scr{B}(H,\breve{L})$ agree. Consequently, by the criterion for special points in the split case
		\ccite{kaletha_prasad}{Lemma~6.4.2}, it suffices to demonstrate that there exists a finite Galois extension $L/F$ such that for each root $\alpha$, with respect to a fixed maximal torus, there is an affine functional with derivative $\alpha$ that vanishes at the vertex $x\in\scr{B}(H,\breve{L})$. This can be done by following the arguments of \ccite{larsen_maximality}{Lemma~2.4}, where loc.~cit.~shows that one may choose any totally ramified, finite Galois extension $L/F$ so that $[L:F]$ is divisible by (equivalently, divides a common multiple of) the denominators occurring in the expression of $x$ as a $\mathbb{Q}$-linear combination of coroots. In particular, when $H$ is semisimple with no bad simple factor, the extension $L/K$ can be chosen to be tamely ramified.
	\es
	The following was proved in \ccite{pappas-rapoport_parahoric_grothendieck_serre}{Proposition~2.7} in the simply connected semisimple case. Our proof is essentially the same as that in loc.~cit.
	\bp\label{prop:larsen_hyperspeciality}
		Let $K$ be a henselian discretely valued field with a perfect residue field. Given a reductive $K$-group scheme $H$ and a $\Gal(\breve{K}/K)$-invariant parabolic subgroup $\ca{P}\subseteq H(\breve{K})$ in the sense of Definition~\ref{defn:parahoric}\eqref{defn:parahoric_subgroup}, there exist 
		\bun 
			\item[$\circ$] a finite Galois extension $L/K$ of discretely valued fields splitting $H$, and
			\item[$\circ$] a reductive $\ca{O}_L$-model $\scr{H}$ of $H\times_{\spec(K)}\spec(L)$ 
		\eun so that we have \bud\ca{P}=\scr{H}(\ca{O}_{\breve{L}})\cap H(\breve{K}).\eud
		Moreover, if $H$ is semisimple and has no bad simple factor in the sense of Definition~\ref{defn:Serre_number}, there is a tamely ramified extension $L/K$ satisfying the claim.
		Furthermore, any such $\ca{O}_L$-group scheme $\scr{H}$ comes equipped with a $\Gal(L/K)$-action compatible with the natural semi-linear $\Gal(L/K)$-action on its generic fibre.
	\ep
	\bs
		We begin by showing the first two claims. By definition, we choose a point $x\in\scr{B}(H,\breve{K})$ such that \bd\label{diag:proof_parahoric}\ca{P}=\text{Stab}_{H(\breve{K})}(x)\cap H(\breve{K})^0.\ed By Proposition~\ref{lem:larsen_hyperspecial}, there exists a finite Galois extension $L/K$ splitting $x$ such that $x$ maps to a hyperspecial point in $\scr{B}(H,\breve{L})$. Moreover, in the case where $H$ is semisimple and has no bad simple factor, this extension $L/K$ may be chosen to be tamely ramified.
		\par Therefore, we note that by hyperspeciality of $x$ and \ccite{kaletha_prasad}{Proposition~8.4.14}, there exists a reductive $\ca{O}_L$-model $\scr{H}$ of $H\times_{\spec(K)}\spec(L)$ such that \bd\label{diag:proof_hyperspecial}\scr{H}(\ca{O}_{\breve{L}})=\Stab_{H(\breve{L})}(x)\cap H(\breve{L})^0.\ed We note that, since $H$ splits over $L$, loc.~cit.~ensures that the reductive $\ca{O}_L$-group scheme $\scr{H}$ is fibrewise connected. Therefore, it remains to verify the displayed equality in the first claim. By Remark~\ref{rem:functoriality_of_Kottwitz}, since we have $$H(\breve{K})^0= H(\breve{K})\cap H(\breve{L})^0,$$
		it suffices to combine the equalities \eqref{diag:proof_parahoric} and \eqref{diag:proof_hyperspecial} with \eqref{diag:building_extension_base_field}. Indeed, \begin{align*}\scr{H}(\ca{O}_{\breve{L}})\cap H(\breve{K}) & =(\mathrm{Stab}_{H(\breve{L})}(x)\cap H(\breve{L})^0)\cap H(\breve{K})\\ &=(\mathrm{Stab}_{H(\breve{K})}(x)\cap H(\breve{K}))\cap H(\breve{K})^0.\end{align*} 
		\par Moreover, since $x\in\scr{B}(H,\breve{L})$ is a $\Gal(L/K)$-invariant point, by \eqref{diag:proof_hyperspecial}, it follows that $\scr{H}(\ca{O}_{\breve{L}})$ is $\Gal(L/K)$-invariant. Consequently, the last claim follows from \ccite{kaletha_prasad}{Corollary~2.10.10}.
	\es
	The following theorem builds on previous known cases in the literature (see Remark~\ref{rem:cor:known_cases}). As mentioned in the introduction, its novelty is the fact that it allows for possibly wildly ramified base field extensions.
	\bt\label{thm:reductification_parahoric}
		Given a parahoric group scheme $\scr{P}$ over a henselian discretely valued field $K$ with a perfect residue field, there exist
		\bun
			\item[$\circ$] a finite Galois extension $L/K$ with Galois group $\mathit{\Gamma}$ and
			\item[$\circ$] a reductive $\ca{O}_L$-model $\scr{G}$ of $\scr{P}\times_{\spec(\ca{O}_K)}\spec(L)$ that respects the natural semilinear $\mathit{\Gamma}$-action on its generic fibre,
		\eun
		such that we have an isomorphism of $\ca{O}_K$-group schemes 
		\bud\scr{P}\cong \fr{R}_{\ca{O}_L/\ca{O}_K}^{\mathit{\Gamma}}(\scr{G})_{\sm}. \eud 
	\et
	\bs
		Thanks to Proposition~\ref{prop:larsen_hyperspeciality}, there exist a finite Galois extension $L$ of $K$ with Galois group $\Gamma$ and a $\Gamma$-equivariant, reductive $\ca{O}_L$-model $\scr{G}$ of $\scr{P}\times_{\spec(\ca{O)}}\spec(L)$ such that \bd\label{diag:them:proof}\scr{P}(\ca{O}_{\breve{K}})=\scr{G}(\ca{O}_{\breve{L}})\cap\scr{P}(\breve{K}).\ed
		Whereas, by \ccite{conrad_gabber_prasad_pseudo_reductive_groups}{Propostion~A.5.2}, the Weil restriction $\fr{R}_{\ca{O}_L/\ca{O}_K}(\scr{G})$ is smooth, affine and is compatible with base change. Therefore, letting $\scr{G}'\colonequals \fr{R}_{\ca{O}_L/\ca{O}_K}^{\mathit{\Gamma}}(\scr{G})$, it suffices to show that $\scr{P}\cong \scr{G}'_{\sm}$. By \ccite{kaletha_prasad}{Corollary~2.10.11}, to prove this claim, it suffices to show that $\scr{P}(\ca{O}_{\breve{K}})\cong\scr{G}'_{\sm}(\ca{O}_{\breve{K}})$. 
		However, since Proposition~\ref{prop:smoothening_exists} ensures that $\scr{G}'_{\sm}(\ca{O}_{\breve{K}})\cong\scr{G}'(\ca{O}_{\breve{K}})$, we reduce to establishing the equality \bd\label{diag:proof:thm:need_to_show}\scr{P}(\ca{O}_{\breve{K}})\cong\scr{G}'(\ca{O}_{\breve{K}}).\ed 
		Whereas, by definition, we have $\scr{G}'(\ca{O}_{\breve{K}})=(\scr{G}(\ca{O}_{\breve{K}}\otimes_{\ca{O}_K}\ca{O}_L))^\Gamma$. 
		Let $K\subseteq K^{u}\subseteq L$ be the maximal unramified subextension and let $I=\Gal(L/K^u)$ and $\Xi=\Gal(K^u/K)$. Then, there is an exact sequence of Galois groups \bd\label{diag:proof:thm:Inertia_subgroup}1\tir&I\tir& \Gamma\tir&\Xi\tir& 1.\ed 
		However, supposing that $\pi\in L$ is a uniformiser, since $L/K^u$ is totally ramified, we can write $\ca{O}_{K^u}[\pi]=\ca{O}_L$, whence, there is an equality \bd\label{diag:proof:thm:unramified_sh}\ca{O}_{\breve{K}}[\pi]=\ca{O}_{\breve{L}}.\ed In particular, by definition, we have that $\ca{O}_{\breve{K}}=(\ca{O}_{\breve{L}})^I$.
		Therefore, there is a $\Gamma$-equivariant isomorphism \begin{align*}\displaystyle\ca{O}_{\breve{K}}\otimes_{\ca{O}_K}\ca{O}_L&\cong \left(\prod_{\Xi}\ca{O}_{\breve{K}}\right)\otimes_{\ca{O}_{K^u}}\ca{O}_L\\&\cong\left(\prod_{\Xi}\ca{O}_{\breve{K}}\right)[\pi]\\&\cong\left(\prod_{\Xi}(\ca{O}_{\breve{K}}[\pi])\right),\end{align*}where the product is equipped with the diagonal $\Xi$-action. We note the $I$-action induced by \eqref{diag:proof:thm:Inertia_subgroup} does not commute the factors of the product above.
		In particular, taking $I$-invariants, we have a $\Xi$-equivariant isomorphism \bud(\scr{G}(\ca{O}_{\breve{K}}\otimes_{\ca{O}_K}\ca{O}_L))^{I}\cong\prod_{\Xi}\scr{G}(\ca{O}_{\breve{K}}[\pi])^I\underset{\eqref{diag:proof:thm:unramified_sh}}{\cong}\prod_{\Xi}\scr{G}(\ca{O}_{\breve{K}}). \eud
		Taking $\Xi$-invariants in the above displayed isomorphism, we obtain \bud(\scr{G}(\ca{O}_{\breve{K}}\otimes_{\ca{O}_K}\ca{O}_L))^{\Gamma}\cong\left(\prod_{\Xi}\scr{G}(\ca{O}_{\breve{K}})\right)^{\Xi}\cong\scr{G}(\ca{O}_{\breve{K}}).\eud This finishes the proof.
	\es
	\br[Fibrewise connected]
		By the uniqueness of integral models proved in \ccite{kaletha_prasad}{Corollary~2.10.11}, we note that the smooth invariant pushforward $\fr{R}_{\ca{O}_L/\ca{O}_K}^{\mathit{\Gamma}}(\scr{G})_{\sm}$ is automatically fibrewise connected.
	\er
	\br[Perfect residue field]
		As remarked in Remark~\ref{rem:perfect_residue_field1}, the approach of Kaletha--Prasad allows us to remove perfect residue field condition provided that $G$ is quasi-split over $\breve{K}$ (see \ccite{kaletha_prasad}{Introduction, Page~14}).
	\er
	We conclude this section with several examples and remarks. In the case of tori, reductifications are closely related to their Néron models. This situation is discussed in Remark~\ref{rem:torus_reductification}, which shows that, in the wildly ramified case, group smoothening is essential (see Example~\ref{exmp:edixhoven}). Finally, Example~\ref{exmp:iwahori} examines reductifications of Iwahori subgroups of $\SL_2$; whereas, Remark~\ref{rem:iwahori_not_tame} shows that base change does not preserve the property of being parahoric.
	\br[Torus case]\label{rem:torus_reductification}
	In the case when $G=T$ is a $K$-torus, the statement of Theorem~\ref{thm:reductification_parahoric} can be made more explicit. Indeed, in this case there is a unique parahoric $\ca{O}_K$-model of $T$ given by the identity component $(T^{\mathrm{ft}})^{\circ}$ of the \textit{ft} Néron model $T^{\mathrm{ft}}$ of $T$ (see \ccite{kaletha_prasad}{\textsection B.4,\textsection B.7 and \textsection B.8} or \ccite{chai_yu_neron_model_tori}{\textsection3}). Letting $L/K$ be a splitting field of $T$, this $T^{\mathrm{ft}}$ can be written as the group smoothening of the schematic closure of $T\subseteq\fr{R}_{\ca{O}_L/\ca{O}_K}(\bb{G}_{m,\ca{O}_L})$. We claim that the canonial morphism \bud (T^{\mathrm{ft}})^{\circ}\cong  \fr{R}^{\mathit{\Gamma}}_{\ca{O}_L/\ca{O}_K}(\bb{G}_{m,\ca{O}_L})_{\sm}\eud induces an isomorphism of smooth $\ca{O}_K$-group schemes with generic fibre $T(K)$. By definition, since the source is smooth over $\ca{O}_K$, thanks to \ccite{kaletha_prasad}{Corollary~2.10.11}, this isomorphism follows by looking at the $\ca{O}_K$-sections.
	\er
	\begin{exmp}[\citeauthor{edixhoven_neron_model_tame_ramification}'s example]\label{exmp:edixhoven}
		For this example, we set $K=\bb{Q}_2[i]$ to be the splitting field of $t^2+1$ over $\bb{Q}_2$. Let $T$ be the twisted form of $\bb{G}_{m,\bb{Q}_2}=\spec(\bb{Q}_2[X^{\pm 1}])$ which is defined to be $\fr{R}^{\mathit{\Gamma}}_{K/\bb{Q}_2}(\bb{G}_{m,K})$, where the action of the cyclic group $\Gamma=\langle\sigma\rangle$ of order $2$ is given by $\sigma\colon i\mapsto -i$ and $\sigma\colon X\mapsto X^{-1}$. By Remark~\ref{rem:torus_reductification}, the parahoric model $\scr{T}$ of $T$ is given by $\fr{R}^{\mathit{\Gamma}}_{\ca{O}_K/\bb{Z}_2}(\bb{G}_{m,\ca{O}_K})_{\sm}$. However, as noted in \ccite{edixhoven_neron_model_tame_ramification}{Example~4.3}, the group smoothening is necessary because $$\fr{R}^{\mathit{\Gamma}}_{\ca{O}_K/\bb{Z}_2}(\bb{G}_{m,\ca{O}_K})=\spec\left(\bb{Z}_2[X,Y]/(X^2+Y^2-1)\right)~~~~~~~~\text{whose special fibre is isomorpic to}~~~~~\bb{G}_{a,\bb{F}_2}[z]/(z^2).$$
	\end{exmp}
	\begin{exmp}[Iwahori subgroup]\label{exmp:iwahori}
		Let $p$ be an arbitrary prime number and let $\Gamma$ be the Galois group of the quadratic extension $K\colonequals \bb{Q}_p(\sqrt{p})/\bb{Q}_p$ with ring of integers $\bb{\ca{O}}\subset K$. The Iwahori subgroup $\ca{I}\subset\SL_2(\bb{Q}_p)$ is given by the set of matrices \[
		\ca{I}\colonequals\left\{\begin{pmatrix}
			\ast & \ast \\
			p\ast & \ast
		\end{pmatrix}\right\}\subset\SL_2(\bb{Z}_p)
		\] whose reduction modulo $p$ is the Borel subgroup given by the upper triangular matrices. Then we can check that 
		$\ca{I}=\fr{R}_{\ca{O}/\bb{Z}_p}^{\mathit{\Gamma}}(\scr{H})(\bb{Z}_p)$, where $\scr{H}$ is the reductive $\ca{O}$-model of $SL_{2,K}$, conjugate to $SL_{2,\ca{O}}$ (since every point is hyperspecial), whose integers points are given by \[
		\scr{H}(\ca{O})=\left\{\begin{pmatrix}
			\ca{O} & \frac{1}{\sqrt{p}}\ca{O} \\
			\sqrt{p}\ca{O} & \ca{O}
		\end{pmatrix}\right\}\subset\SL_2(K)\]
		Thus $(K,\scr{H})$ is a reductification and since $K/\bb{Q}_p$ is quadratic extension, it is even a minimal reductification of $\ca{I}$.
		\par In the case when $p=2$, the extension $K/\bb{Q}_2$ is wildly ramified. In this situation, it is natural to ask the following questions, whose answer is unclear to the author. Is the group scheme $\fr{R}_{\ca{O}/\bb{Z}_2}^{\mathit{\Gamma}}(\scr{H})$ automatically smooth? Does $\ca{I}$ have any tame reductification?
	\end{exmp}
	\br[Base change of parahoric]\label{rem:iwahori_not_tame}
		Let $L/K$ be a finite extension of discretely valued fields and $\scr{P}$ be a parahoric group scheme over $\ca{O}_L$. It is true in the case when $L/K$ is unramified that $\scr{P}\times_{\spec(\ca{O}_K)}\ca{O}_L$ remains a parahoric group. However, this is well known that this is no longer true even when $L/K$ is tamely ramified. For example, we take $K=\bb{Z}_p$ and a finite extension $L$ with uniformiser $\varpi$. Then, as remarked in Example~\ref{exmp:iwahori}, the Iwahori subgroup is the non-reductive parahoric model $\scr{I}_{L}\subset\SL_{2,\ca{O}_L}$ whose set of integers points is given by the set of matrices \[
		\scr{I}_L(\ca{O}_L)\colonequals\left\{\begin{pmatrix}
			\ast & \ast \\
			\varpi\ast & \ast
		\end{pmatrix}\right\}\subset\SL_2(\ca{O}_L).
		\]
		As long as $\varpi\neq p$, equivalently, as long as the extension $L/K$ is ramified,\footnote{The author thanks Paul Broussous for highlighting this example in MathOverflow post \href{https://mathoverflow.net/questions/459156/arbitrary-base-change-of-a-parahoric-subgroup-in-split-case}{\#459156}.} the above set is clearly different from \bud(\scr{I}_{K}\times_{\spec(\ca{O}_K)}\spec(\ca{O}_L))(\ca{O}_L)=\scr{I}_K(\ca{O}_L).\eud
	\er
	\section{Tamely reductifiable parahoric group schemes.}\label{sec:gs_serre}
	In this section, our goal is to prove Corollary~\ref{cor:general_case_gs}, which establishes the Grothendieck--Serre conjecture for parahoric group schemes which are generically simply-connected and in good residue characteristics. This result is a consequence of two statements. First, in Theorem~\ref{thm:gs_parahoric}, we demonstrate the triviality of any generically trivial torsor under any tamely reductifiable parahoric group scheme. Our strategy follows the steps of the proof in the reductive case in \cite{ning_dedekind}. Second, in Corollary~\ref{cor:tame_reductification_parahoric}, we show that a parahoric group scheme is tamely reductifiable in good residue characteristics, thereby establishing the second claim of Theorem~\ref{thm:B}. This extends various results from the literature (see Remark~\ref{rem:cor:known_cases}). We obtain Corollary~\ref{cor:tame_reductification_parahoric} as a consequence of Theorem~\ref{thm:reductification_parahoric} and Proposition~\ref{prop:tamely_reductifiable_z-extension}.
	\par We begin with the following definition and lemma, which are required to state Proposition~\ref{prop:tamely_reductifiable_z-extension}.
	\bdf[Centrally tame]\label{defn:centrally_tame}
	A torus $T$ over a discretely valued field $K$ is called \textit{tame} if its splits over a tamely ramified extension of $K$ and a reductive $K$-group scheme $G$ is called \textit{centrally tame} if its maximal central torus is tame.
	\edf
	\br\label{rem:tame_torus}
		A tame torus is tamely reductifiable. Indeed, this follows from Remark~\ref{rem:torus_reductification}.
	\er
	\bl\label{lem:properties_tame_reductification}
		Let $K$ be a henselian discretely valued field and let $\scr{P}$ be a tamely reductifiable parahoric $\ca{O}_K$-group scheme. Then, the base change $(F,\scr{G}_{\ca{O}_F})$ of any tame reductification $(L,\scr{G})$ of $\scr{P}$ along a tamely ramified extension $F/L$ remains a tame reductification of $\scr{P}$. In particular, given another tamely reductifiable parahoric $\ca{O}_K$-group scheme $\scr{P}'$, the product $\scr{P}\times\scr{P}'$ is tamely reductifiable.
	\el
	\bs
		We start with the proof of the first claim. By \ccite{kaletha_prasad}{Corollary~2.10.11}, we note that the canonical morphism $\scr{G}\to\fr{R}^{\Gal(F/L)}_{\ca{O}_F/\ca{O}_L}(\scr{G}_{\ca{O}_F})$ is an isomorphism. Consequently, the claim follows if we establish that the canonical morphism \bud\fr{R}^{\Gal(F/K)}_{\ca{O}_F/\ca{O}_K}(\scr{G}_{\ca{O}_F})\cong\fr{R}^{\Gal(L/K)}_{\ca{O}_L/\ca{O}_K}(\fr{R}^{\Gal(F/L)}_{\ca{O}_F/\ca{O}_L}(\scr{G}_{\ca{O}_F}))\eud is an isomorphism. In a similar vein as the above, this is a consequence of loc.~cit.  
		\par Thanks to the first claim, to prove the second claim, without loss of generality, it suffices to assume that both $\scr{P}$ and $\scr{P}'$ are tamely reductifiable over $L$. Therefore, this claim can be proven by noting that the functor $\fr{R}^{\Gal(L/K)}_{\ca{O}_L/\ca{O}_K}(-)$ commutes with products (see \ccite{conrad_gabber_prasad_pseudo_reductive_groups}{Proposition~A.5.2(3)} and \ccite{edixhoven_neron_model_tame_ramification}{Proposition~3.1}).
	\es
	\bp\label{prop:tamely_reductifiable_z-extension}
	Let $K$ be a henselian discretely valued field and a $\scr{P}$ be a parahoric $\ca{O}_K$-group scheme with generic fibre $\scr{P}_K$. Let $x\in\scr{B}(\scr{P}_K,\breve{K})$ be a point associated to $\scr{P}$ (see Definition~\ref{defn:parahoric}). We consider the following statements.
	\bn[(i)] 
	\item\label{prop:tamely_ramified_z-extension:pt:1:} The group $\scr{P}_K$ is centrally tame and the simply-connected cover of its derived subgroup has no bad simple factors in the sense of Definition~\ref{defn:Serre_number}.
	\item\label{prop:tamely_ramified_z-extension:pt:2:} There is a $\mu$-extension $1\to\mu\to\tilde{G}\to\scr{P}_K\to 1$ over $K$, in the sense of Definition~\ref{defn:z-extension}, such that $\tilde{G}$ has a tamely reductifiable parahoric $\ca{O}_K$-model $\tilde{\scr{G}}$ which is associated the point in the preimage of $x$ under \eqref{diag:building_quotient_map}.
	\item\label{prop:tamely_ramified_z-extension:pt:3:} This $\scr{P}$ is tamely reductifiable. 
	\en 
	Then \eqref{prop:tamely_ramified_z-extension:pt:1:}$\implies$\eqref{prop:tamely_ramified_z-extension:pt:2:}$\implies$\eqref{prop:tamely_ramified_z-extension:pt:3:}.
	\ep
	\bs
	\eqref{prop:tamely_ramified_z-extension:pt:1:}$\implies$\eqref{prop:tamely_ramified_z-extension:pt:2:}: Let $G_{\SC}$ be the simply-connected cover of the derived subgroup of $\scr{P}_K$. Using the identification $\scr{B}(G_{\SC},\breve{K})=\scr{B}(\scr{P}_K,\breve{K})$, proved in \ccite{kaletha_prasad}{(7.6.2.1)}, let $\scr{G}$ be the parahoric model of $G_{\SC}$ corresponding to $\scr{P}$. By assumption, the semisimple group $G_{\SC}$ has no bad simple factors; in particular, this $\scr{G}$ is tamely reductifiable. Let \bud1\to\mu\to\tilde{G}\to\scr{P}_K\to 1\eud be the $\mu$-extension constructed in \eqref{rem:mu_extension}, i.e., $\tilde{G}=G_{\SC}\times T$, where $T$ is the maximal central torus of $\scr{P}_K$. Since $T$ splits over a tamely ramified extension of $K$, the parahoric model $\scr{T}$ of $T$ is tamely reductifiable (see Remark~\ref{rem:tame_torus}). Consequently, by Lemma~\ref{lem:properties_tame_reductification}, the parahoric group scheme $\tilde{\scr{G}}\colonequals\scr{G}\times\scr{T}$ is tamely reductifiable. Finally, to show that $\tilde{\scr{G}}$ is associated to the point in the preimage of $x$, we employ the equality $\scr{B}(\tilde{G},\breve{K})=\scr{B}(\scr{P}_K,\breve{K})$, the existence of which follows from \eqref{diag:building_quotient_map} and \ccite{kaletha_prasad}{(7.6.2.1)}.
	\vspace{0.15cm}\par\eqref{prop:tamely_ramified_z-extension:pt:2:}$\implies$\eqref{prop:tamely_ramified_z-extension:pt:3:}: Let $\tilde{\scr{G}}$ be a tamely reductifiable parahoric $\ca{O}_K$-model of $\tilde{G}$ associated to the point in the preimage of $x$, let $(K',\tilde{\scr{G}}')$ be a tame reductification of $\tilde{\scr{G}}$ and let $\Gamma\colonequals\Gal(K'/K)$. By Remark~\ref{rem:z-extension}\eqref{rem:z-extension:pt:1}, we note that \bud\text{the Zariski-closure $\tilde{\mu}'$ of $\mu\times_{\spec(K)}\spec(K')$ inside $\tilde{\scr{G}}'$}\eud lies inside the centre $Z(\tilde{\scr{G}}')$ of the reductive $\ca{O}_{K'}$-group scheme $\tilde{\scr{G}}'$. Indeed, this follows from the flatness of $Z(\tilde{\scr{G}}')$ over $\ca{O}_{K'}$. We equip $\tilde{\mu}'$ with the canonical $\Gamma$-semilinear action. By \ccite{conrad_reductive_group_schemes}{Theorem~3.3.4}, this $Z(\tilde{\scr{G}}')$ is of multiplicative type; consequently, by \ccite{conrad_reductive_group_schemes}{Corollary~B.3.3}, the same holds for $\tilde{\mu}'$. Therefore, the quotient \bud\scr{P}'\colonequals \tilde{\scr{G}}'/\tilde{\mu}'\eud is a reductive $\ca{O}_K$-group scheme (see, for example, \ccite{sga3ii}{exposé~IX, proposition~2.3}), which borrows a canonical $\Gamma$-equivariant structure. Taking invariant pushforwards, we obtain a natural triangle of morphisms of smooth $\ca{O}_K$-group schemes (for example, by looking at the $\ca{O}_{\breve{K'}}$-sections and \ccite{kaletha_prasad}{Lemma~2.10.13}) \bud\tilde{\scr{G}}\arrow[rr, two heads]\arrow[dr]&&\scr{P}\\&\fr{R}^{\mathit{\Gamma}}_{\ca{O}_{K'}/\ca{O}_K}(\scr{P}')\arrow[ur,swap,"\varphi"].\eud To prove our claim, it suffices to show that $\varphi$ is an isomorphism. By \ccite{kaletha_prasad}{Corollary~2.10.11} and the smoothness of both group schemes, this follows from the bijection at the level of $\ca{O}_{\breve{K}}$-points.
	\es
	\bc\label{cor:tame_reductification_parahoric}
	Given a henselian discretely valued field $K$ and a parahoric $\ca{O}_K$-group scheme $\scr{P}$ 
	\bn[(a)] 
		\item\label{cor:pt:a} whose generic fibre $\scr{P}_K$ is either an adjoint semisimple $K$-group scheme that has no bad simple factor, or 
		\item\label{cor:pt:b} such that $\scr{P}_K$ is centrally tame and the simply-connected cover of its derived subgroup has no bad simple factor,
	\en 
	there exist
	\bun
	\item[$\circ$] a tamely ramified finite Galois extension $L/K$ with Galois group $\mathit{\Gamma}$ and
	\item[$\circ$] a reductive $\ca{O}_L$-model $\scr{G}$ of $\scr{P}\times_{\spec(\ca{O}_K)}\spec(L)$ that respects the natural semilinear $\mathit{\Gamma}$-action on its generic fibre,
	\eun
	such that we have an isomorphism of $\ca{O}_K$-group schemes 
	\bud\scr{P}\cong \fr{R}_{\ca{O}_L/\ca{O}_K}^{\mathit{\Gamma}}(\scr{G}). \eud 
	\ec
	\bs
		Case \eqref{cor:pt:a} follows from Theorem~\ref{thm:reductification_parahoric} and the result \ccite{edixhoven_neron_model_tame_ramification}{Proposition~3.4}, which states that when $L/K$ is tamely reductifiable, the $\ca{O}_K$-scheme $\fr{R}_{\ca{O}_L/\ca{O}_K}^{\mathit{\Gamma}}(\scr{G})$ is already smooth. Case \eqref{cor:pt:b} follows from \eqref{cor:pt:a} and Proposition~\ref{prop:tamely_reductifiable_z-extension}.
	\es
	\br\label{rem:cor:known_cases}
	Corollary~\ref{cor:tame_reductification_parahoric} builds on prior results in the literature. Balaji and Seshadri proved the statement in residue characteristic~$0$ in \ccite{balaji_seshadri}{Theorem~5.2.7}, under the assumption that the generic fibre $G_K\colonequals\scr{P}\times_{\spec(\ca{O}_K)}\spec(K)$ is a base change of a reductive group defined over $\bb{Z}$. Generalisations of this results were provided in \cite{damioli_thesis, damiolini_hong_local_types_equivariant_bundles}. In positive and mixed characteristics, Pappas and Rapoport established the result \ccite{pappas-rapoport_parahoric_grothendieck_serre}{Proposition~2.8} for the case where $G_K$ is a simply-connected, absolutely simple, semisimple group scheme (see \cite{chernousov_gille_pianzola_parahoric} for a related result).
	\er
		We now turn to the proof of Theorem~\ref{thm:gs_parahoric}, which is a consequence of Propositions~\ref{prop:ning}, Lemma~\ref{lem:group_to_parabolic} and Proposition~\ref{lem:parabolic_to_Levi} stated below.
	
	\bp\label{prop:ning}
		Given a henselian discrete valuation field $K$ with a perfect residue field and a parahoric $\ca{O}_K$-group scheme $\scr{P}$ whose generic fibre $\scr{P}_K$ is anisotropic and simply-connected semisimple, any generically trivial $\scr{P}$-torsor over $\ca{O}_K$ is trivial. 
	\ep
	\bs
		We closely follow the proof in \ccite{ning_dedekind}{Theorem~5.1}. Let $\breve{K}$ be a maximal unramified extension contained in a separable closure $K^{\sep}$ of $K$. We suppose that $\Gamma\colonequals\Gal(K^{\sep}/K)$ and $\Xi\colonequals\Gal(\breve{K}/K)$. To prove the claim, it suffices to show that the restriction morphism \bud H^1(\ca{O}_K,\scr{P})\to H^1(K,\scr{P}_K)\eud has trivial kernel. Since $K$ is a field, the target cohomology group is isomorphic to $H^1(\Gamma,\scr{P}(K^{\sep}))$. Furthermore, since $\ca{O}_K$ is henselian, we have $H^1(\ca{O}_K,\scr{P})\cong H^1(\Xi,\scr{P}(\ca{O}_{\breve{K}}))$. Thus, it suffices to show the equality \bud \ker(H^1(\Xi,\scr{P}(\ca{O}_{\breve{K}}))\to H^1(\Gamma,\scr{P}(K^{\sep})))=\{\ast\}.\eud
		By the inflation-restriction exact sequence \ccite{serre_galois_cohomology}{Chapter~I, \textsection5.8 a)}, the canonical morphism $H^1(\Xi,\scr{P}(\breve{K}))\to H^1(\Gamma,\scr{P}(K^{\sep}))$ has trivial kernel. Therefore, it suffices to show the following equality \bud \ker(H^1(\Xi,\scr{P}(\ca{O}_{\breve{K}}))\to H^1(\Xi,\scr{P}(\breve{K})))=\{\ast\}.\eud Let $c\in H^1(\Xi,\scr{P}(\ca{O}_{\breve{K}}))$ be an arbitrary $1$-cocyle that maps to a coboundary in $H^1(\Xi,\scr{P}(\breve{K}))$, i.e., there exists an element $g\in\scr{P}(\breve{K})$ such that $c(\gamma)=g^{-1}(\gamma\cdot g)$ for all $\gamma\in\Xi$. It suffices to show that $c$ itself is a coboundary; in particular, it is enough to establish that $g\in\scr{P}(\ca{O}_{\breve{K}})$, which we show below. \par We note that the simply-connectedness of $\scr{P}_K$ implies that $\scr{P}(\breve{K})^0=\scr{P}(\breve{K})$ (see \eqref{defn:G(K)^0}). As a consequence, by Definition~\ref{defn:parahoric}\eqref{defn:parahoric:1}, there is a point $x\in\scr{B}(\scr{P}_K,\breve{K})$ such that $\scr{P}(\ca{O}_{\breve{K}})$ is the stabiliser of $x$ in the $\scr{P}(\breve{K})$-action on $\scr{B}(\scr{P}_K,\breve{K})$. Consequently, since $\scr{P}$ is defined over $\ca{O}_K$, this $x$ is $\Xi$-invariant. Let us consider the point $g\cdot x\in\scr{B}(\scr{P}_K,\breve{K})$, whose stabiliser is $g\scr{P}(\ca{O}_{\breve{K}})g^{-1}$. By the definition of $g$, a computation shows that $g\cdot x$ is also $\Xi$-invariant. However, since $\scr{P}$ is anisotropic, thanks to \ccite{kaletha_prasad}{Proposition~9.3.9}, there is a unique $\Gamma$-invariant point in $\scr{B}(\scr{P}_K,\breve{K})$, proving that $x=g\cdot x$. In particular, we have that \bud\scr{P}(\ca{O}_{\breve{K}})=g\scr{P}(\ca{O}_{\breve{K}})g^{-1}.\eud Equivalently, the above displayed equality shows that $g$ normalises $\scr{P}(\ca{O}_{\breve{K}})$. However, since $\scr{P}_K$ is simply-connected, by \ccite{kaletha_prasad}{\textsection9.3.26}, we know that $\scr{P}(\ca{O}_{\breve{K}})$ is its own normaliser in $\scr{P}(\breve{K})$. Thus, we obtain that $g\in\scr{P}(\ca{O}_{\breve{K}})$, and the proof is complete.
	\es
	\bl\label{lem:group_to_parabolic}
		Given a henselian discretely valued field $L$, a profinite group $\Gamma$ which is a cofiltered limit of finite subgroups of the group of automorphisms of $\ca{O}_L$, and a reductive $\ca{O}_L$-group scheme $\scr{G}$ equipped with a semilinear $\Gamma$-action, every $\Gamma$-equivariant parabolic subgroup $\scr{Q}_L\subset\scr{G}_L$ of the generic fibre lifts to a $\Gamma$-equivariant parabolic subgroup $\scr{Q}\subset\scr{G}$, respecting the $\Gamma$-action on the generic fibre. Furthermore, any such lift induces an equivalence between the categories of generically trivial $\Gamma$-equivariant torsors under $\scr{Q}$ and that under $\scr{G}$.
	\el
	\bs
		By the valuative criterion of properness applied to the scheme parametrising proper parabolic subgroups of $\scr{G}$, any such parabolic subgroup $\scr{Q}_L\subsetneq\scr{G}_L$ uniquely lifts to a parabolic subgroup $\scr{Q}\subsetneq\scr{G}$. Furthermore, there is even an equality of points $(\scr{G}/\scr{Q})(\ca{O}_L)=(\scr{G}/\scr{Q})(L)$, which, in turn, proves the equality $\scr{Q}(\ca{O}_L)=\scr{Q}(L)\cap\scr{G}(\ca{O}_L)$. Consequently, by \ccite{kaletha_prasad}{Corollary~2.10.10}, this $\scr{Q}$ can be endowed with a unique $\Gamma$-action compatible with that on $\scr{Q}_L$ and $\scr{G}$. The final claim then follows as in the proof of \ccite{phd-thesis}{Proposition~7.1.1}).
	\es
	The next result may be viewed as an equivariant analogue of Mostow's theorem \cite{mostow_levi_decomposition} (which has since been generalised in \cite{mcninch_tame_descent_levi, mcninch_charp_perfect_residue_field_levi, demarche_mostow_levi_decomposition}). Before stating it, let us recall some necessary terminology (cf.~\ccite{sga3iii}{exposé~XXVI, proposition~1.6}).
	\par Let $\scr{Q}$ be a fibrewise connected, smooth, affine group scheme over a scheme $S$, and let $\scr{U}\subseteq \scr{Q}$ denote its unipotent radical. A Levi decomposition of $\scr{Q}$ is the datum of an $S$-subgroup scheme $\scr{M}\subseteq \scr{Q}$ together with a split exact sequence of $S$-group schemes \bud1\tir&\scr{U}\tir&\scr{Q}\tir&\scr{M}\tir\arrow[l, bend right, hook']&1.\eud The subgroup $\scr{M}$ is then called a \textit{Levi subgroup} of $\scr{Q}$. Let $\ca{E}$ be the presheaf on $S$ sending an $S$-scheme $T$ to the set $\ca{E}(T)$ of Levi decomposition of $\scr{Q}$ over $T$. For a $\Gamma$-action on $S$, we write $H^p_{\Gamma}(S,\scr{U})$ for the $p$-th $\Gamma$-equivariant étale cohomology group of $S$ with coefficients in $\scr{U}$.
	\bp\label{lem:parabolic_to_Levi}
	Let $L$ be a henselian discretely valued field of residue characteristic $p$, let $\Gamma$ be a finite subgroup of the group of automorphisms of $\ca{O}_L$ whose order is coprime to $p$ and let $\scr{Q}$ be a fibrewise connected, smooth, affine $\ca{O}_L$-group scheme $\scr{Q}$ equipped with a semilinear $\Gamma$-action. Assume that 
	\bn[(i)] 
		\item\label{i} the unipotent radical $\scr{U}\subseteq\scr{Q}$ admits a $\Gamma$-equivariant filtration by subgroups \bud 1=\scr{U}_0\subsetneq\scr{U}_1\subsetneq\ldots\subsetneq\scr{U}_m=\scr{U}\eud whose successive quotients $\scr{U}_n/\scr{U}_{n-1}$ are vector groups $\bb{V}(\scr{F}_n)$ on which $\Gamma$-acts semilinearly, and
		\item\label{ii} the presheaf $\ca{E}$ parametrising Levi decompositions of $\scr{Q}$ is a $\scr{U}$-torsor over $\ca{O}_L$.
	\en 
	Then, there exists a $\Gamma$-equivariant Levi decomposition \bd\label{diag:equivariant_levi}1\tir&\scr{U}\tir&\scr{Q}\tir&\scr{M}\tir\arrow[l, bend right]&1.\ed Any such exact sequence as above induces an equivalence between the categories of (resp., generically trivial) $\Gamma$-equivariant torsors under $\scr{Q}$ and that under $\scr{M}$.
	\ep
	\bs
		 We first demonstrate that \bd\label{diag:vanishing_proof} H^1_{\Gamma}(\ca{O}_L,\scr{U})=0.\ed Consider the local-to-global spectral sequence \ccite{sga4ii}{exposé~V, proposition~6.4} $H^p(\ca{O}_L,\scr{H}^q_{\Gamma}(\ca{O}_L,\scr{U}))\Rightarrow H^{p+q}_{\Gamma}(\ca{O}_L,\scr{U})$, where $\scr{H}^q_{\Gamma}(\ca{O}_L,\scr{U})$ denotes the étale sheafification of the presheaf mapping an $\ca{O}_L$-scheme $T$ to $H^q_{\Gamma}(T,\scr{U})$. It therefore suffices to establish that \bd\label{two_vanishing}\text{$\scr{H}^1_{\Gamma}(\ca{O}_L,\scr{U})=0$~~~~~~~and~~~~~$H^1(\ca{O}_L,\scr{U}^{\Gamma})=0$.}\ed Taking the stalk of $\scr{H}^1_{\Gamma}(\ca{O}_L,
		 \scr{U})$ at the closed point, the first vanishing reduces to showing \bd\label{first_vanishing}H^1_{\Gamma}(\ca{O}_{\breve{L}},\scr{U})=0,\ed which we prove below.
		 \par Since the residue characteristic of $\ca{O}_L$ is coprime to the order of $\Gamma$, the constant group scheme associated to $\Gamma$ is linearly reductive. Consequently, the classifying stack $\mathbf{B}\Gamma$ of $\Gamma$-torsors over $\ca{O}_{\breve{L}}$ is a tame Deligne--Mumford stack in the sense of \cite{abramovich_olsson_angelo_tame_stacks} (see, for example, the paragraph before op.~cit.~Theorem~3.2), and therefore, by the Serre vanishing theorem, the higher coherent cohomology groups of $\mathbf{B}\Gamma$ vanish. Consequently, we have \bd\label{displayed_vanishing}H^q(\mathbf{B}\Gamma,\bb{G}_a)\cong H^q_{\mathrm{coh}}(\mathbf{B}\Gamma,\ca{O}_{\mathbf{B}\Gamma})=0\text{, for all }q>0.\ed Indeed, since $\mathbf{B}\Gamma$ has an étale cover by an affine scheme, the first isomorphism follows by an argument similar to \stacks{03p2}. Therefore, by \eqref{displayed_vanishing} and the fact that $\Gamma$ acts semilinearly on each $\bb{V}(\scr{F}_n)$, we have that $H^q(\mathbf{B}\Gamma, \bb{V}(\scr{F}_n))=0$ for each $q>0$. Consequently, writing the associated long exact sequence of cohomology groups, by induction, it follows that \bd\label{diag:tame_condition} H^q(\mathbf{B}\Gamma,\scr{U})=0\text{, for all }q>0.\ed 
		 However, since $\ca{O}_{\breve{L}}$ is strictly henselian, the $\Gamma$-equivariant small étale topos of $\ca{O}_{\breve{L}}$ is equivalent to the small étale topos of $\mathbf{B}\Gamma$. In particular, there is an isomorphism, $ H^1(\mathbf{B}\Gamma,\scr{U})\cong H^1_{\Gamma}(\ca{O}_{\breve{L}},\scr{U})$, establishing the required vanishing \eqref{first_vanishing}.
		 \par We now turn to prove the second vanishing in \eqref{two_vanishing}. Each $\bb{V}(\scr{F}_n)$ is represented on the étale site by a coherent $\ca{O}_L$-module, hence the same holds for $\scr{U}_n$. Taking $\Gamma$-invariants preserves coherence because $\Gamma$ is finite. Thus $\scr{U}^{\Gamma}$ admits a finite filtration with coherent successive quotients. Consequently, by the Serre vanishing theorem, we obtain $H^1(\ca{O}_L,\scr{U}^{\Gamma})=0$, as required. This shows \eqref{diag:vanishing_proof}. In the rest of the proof, we show that \eqref{diag:vanishing_proof} implies the required claims.
		 \par Since, by assumption, the presheaf $\ca{E}$ is an $\scr{U}$-torsor over $\ca{O}_L$, the presheaf $\ca{E}'$ parametrising $\Gamma$-equivariant Levi decompositions of $\scr{Q}$ is a $(\Gamma,\scr{U})$-bundle in the sense of \ccite{damioli_thesis}{Definition~3}. However, by \eqref{diag:vanishing_proof}, this $\ca{E}'$ has a section over $\ca{O}_L$, showing the existence of \eqref{diag:equivariant_levi}. Finally, suppose that we are given such a Levi decomposition. By the splitness of the exact sequence \eqref{diag:equivariant_levi}, the long exact sequence of $\Gamma$-equivariant cohomology groups yields a short exact sequence \bud 1\to H^1_{\Gamma}(\ca{O}_L,\scr{U})\to H^1_{\Gamma}(\ca{O}_L,\scr{Q})\to H^1_{\Gamma}(\ca{O}_L,\scr{M})\to 1.\eud The final claim then follows from \eqref{diag:vanishing_proof}, completing the proof.
	\es
	As an ingredient in Theorem~\ref{thm:gs_parahoric}, we need to know that an equivariant parabolic subgroup $\scr{Q}$ of a reductive $\ca{O}_L$-group scheme satisfies the hypotheses of Proposition~\ref{lem:parabolic_to_Levi}, which we prove in the following.
	\bl\label{lem:parabolic_splitting}
	Let $L$ be a henselian discretely valued field, let $\Gamma$ be a finite subgroup of the group of automorphisms of $\ca{O}_L$ and let $\scr{Q}$ be a $\Gamma$-equivariant parabolic subgroup $\scr{Q}$ of a reductive $\ca{O}_L$-group scheme endowed with a semilinear $\Gamma$-action. Then, the unipotent radical $\scr{U}\subseteq\scr{Q}$ satisfies conditions \eqref{i} and \eqref{ii} in Proposition~\ref{lem:parabolic_to_Levi}.
	\el
	\bs
	Since $\scr{U}$ is preserved by automorphisms of $\scr{Q}$, it naturally obtains a compatible semilinear $\Gamma$-action.  Furthermore, the canonical filtration on $\scr{U}$ (see \ccite{sga3iii}{exposé~XXVI, proposition~2.1}) also inherits a compatible semilinear $\Gamma$-action. Consequently, this $\Gamma$ acts semilinearly on each graded piece $\bb{V}(\scr{F})$, showing that $\scr{U}$ satisfies \eqref{i}. To show that it satisfies \eqref{ii}, it suffices to use \ccite{sga3iii}{exposé~XXVI, corollaire~1.9}.
	\es
	\bt\label{thm:gs_parahoric}
		Given a henselian discrete valuation ring $\ca{O}_K$ with a perfect residue field and fraction field $K$ as well as a parahoric $\ca{O}_K$-group scheme $\scr{P}$ that is tamely reductifiable in the sense of Definition~\ref{defn:redutification} and whose generic fibre $\scr{P}_K$ is semisimple, there is an equality 
		\bud \ker(H^1(\ca{O}_K,\scr{P})\to H^1(K,\scr{P}_K))=\{\ast\}.\eud 
	\et
	\bs
		Let $(L,\scr{G})$ be a tame reductification of $\scr{P}$ with the Galois group $\Gamma\colonequals\Gal(L/K)$. In this proof, for each $R\in\{\ca{O}_L,L\}$, given a $\Gamma$-equivariant $R$-group scheme $\scr{H}$, the notation $H^1_{\Gamma}(R,\scr{H})$ will denote the pointed set of $\Gamma$-equivariant equivalence classes of $\Gamma$-equivariant $\scr{H}$-torsors over $\spec(R)$.  Since $L/K$ is tame, to prove our claim, thanks to \ccite{pappas-rapoport_parahoric_grothendieck_serre}{Proposition~3.2}, it suffices to demonstrate an equality \bud \ker(H^1_{\Gamma}(\ca{O}_L,\scr{G})\to H^1_{\Gamma}(L,\scr{G}_L))=\{\ast\},\eud where $\scr{G}_L$ denotes the generic fibre of $\scr{G}$. We proceed by inducting on the rank of $\scr{G}_L$.
	\par Thanks to Proposition~\ref{prop:ning}, we can assume that $\scr{P}_K$ is isotropic, which ensures that $\scr{P}_K$ has a proper parabolic subgroup (see \ccite{ces_torsors}{\textsection1.3.6}). By pulling-back a minimal parabolic subgroup of $\scr{P}_K$, we thus obtain a $\Gamma$-equivariant parabolic subgroup $\scr{Q}_L$ of $\scr{G}_L$. By Lemma~\ref{lem:group_to_parabolic}, this subgroup spreads out to a $\Gamma$-equivariant parabolic $\ca{O}_L$-subgroup $\scr{Q}$ of $\scr{G}$. Additionally, the same lemma demonstrates the canonical isomorphism \bud \ker(H^1_{\Gamma}(\ca{O}_L,\scr{Q})\to H^1_{\Gamma}(L,\scr{Q}_L))\cong \ker(H^1_{\Gamma}(\ca{O}_L,\scr{G})\to H^1_{\Gamma}(L,\scr{G}_L)).\eud
	Therefore, to complete the inductive step, it reduces to show the triviality of the kernel on the left.	
		\par By Proposition~\ref{lem:parabolic_to_Levi} and Lemma~\ref{lem:parabolic_splitting}, there is a $\Gamma$-equivariant Levi decomposition $\scr{G}\cong\scr{M}\rtimes\scr{U}$. Furthermore, the same lemma establishes the canonical isomorphism \bud \ker(H^1_{\Gamma}(\ca{O}_L,\scr{Q})\to H^1_{\Gamma}(L,\scr{Q}_L))\cong \ker(H^1_{\Gamma}(\ca{O}_L,\scr{M})\to H^1_{\Gamma}(L,\scr{M}_L)).\eud However, the minimality of the parabolic subgroup $\scr{Q}_L$ guarantees that the rank of the reductive $L$-group $\scr{M}_L$ is strictly smaller than that of $\scr{G}_L$. Consequently, the kernel on the right is trivial by induction. This completes the inductive step, and hence, the proof.
	\es
	Given a central simple algebra $D$ over a henselian discretely valued field $K$ with a residue field $\kappa$, we denote the group of elements of reduced norm one by $\SL_1(D)$ and its adjoint quotient by $\PGL_1(D)$. These groups are non-quasi-split inner $K$-forms of $\SL_m$ and $\PGL_m$, respectively, for a suitable integer $m\geqslant 2$. In particular, when $p\mid m$, these groups are bad simple factors. In fact, if $\cd(\kappa)\leqslant 1$ and $p\neq 2,3,5$, we know that any bad simple factor is isomorphic to a group of this type (see Remark~\ref{rem:delta>1}). The following, which answers Question~\ref{quest:bayer_fluckigar_first} for these groups, is therefore of independent interest.
	\bp\label{prop:A_n-case}
		Let $K$ be a henselian discretely valued field with a perfect residue field $\kappa$ of characteristic $p\neq 2$, let $D$ be a central simple $K$-algebra and let $\scr{P}$ be a parahoric $\ca{O}_K$-group scheme whose generic fibre is isomorphic either to $\SL_1(D)$ or to $\mathrm{PGL}_1(D)$. Then, there exists an $\ca{O}_K$-subgroup scheme $\scr{R}\subseteq\scr{P}$ which is reductive and whose special fibre $\scr{R}_{\kappa}$ canonically maps isomorphically onto the maximal reductive quotient $\scr{P}_{\mathrm{red},\kappa}$ of the $\ca{O}_K$-special fibre $\scr{P}_{\kappa}$ of $\scr{P}$. Furthermore, we have an equality \bd\label{diag:claim_trivial_kernel} \ker(H^1(K,\scr{R})\to H^1(K,\scr{P}))=\{\ast\},\ed
		and in particular, any generically trivial $\scr{P}$-torsor over $\ca{O}_K$ is trivial.
	\ep
	\bs
		Since $\kappa$ is perfect, the maximal $K$-subfield $F$ contained in $D$ is unramified over $K$ (see \ccite{serre_local_fields}{Chapter~XII, \textsection2, Proposition~2}). This $D$ is, in fact, split over $F$, and therefore, so is the generic fibre $\scr{P}_K$ of $\scr{P}$. Consequently, the existence of such an $\scr{R}$ is the content of \ccite{mcninch_reductive_subgroup_parahoric_groups}{Theorem~1}. By \ccite{mcninch_reductive_subgroup_parahoric_groups}{\textsection5.1}, in the case $\scr{P}_K=\SL_1(D)$, the generic fibre $\scr{R}_K$ fits into an exact sequence \bud 1\to\scr{R}_K\to\fr{R}_{F/K}(\bb{G}_{m,F})\to\bb{G}_{m,K}\to 1.\eud Thanks to Hilbert Theorem 90, a calculation then shows that $H^1(K,\scr{R}_K)=\{\ast\}$, whence the equality \eqref{diag:claim_trivial_kernel} follows. Alternatively, when $\scr{P}=\mathrm{PGL}_1(D)$, by definition, there is a morphism $H^1(K,\scr{P}_K)\to\mathrm{Br}(K)$ (with trivial kernel). However, adapting the arguments in loc.~cit., we establish that there is an exact sequence \bud 1\to\bb{G}_{m,K}\to\fr{R}_{F/K}(\bb{G}_{m,F})\to\scr{R}_K\to 1,\eud which shows that there is an injection $H^1(K,\scr{R}_K)\hookrightarrow\mathrm{Br}(K)$, establishing the equality \eqref{diag:claim_trivial_kernel}. It therefore remains to prove the final claim.
		\par By, for example, \ccite{bouthier-ces}{Theorem~2.1.6}, we get an isomorphism $H^1(\ca{O}_K,\scr{P})\cong H^1(\kappa,\scr{P}_{\kappa})$. Since $\kappa$ is perfect, the unipotent radical of $\scr{P}_{\kappa}$ is split, which demonstrates that $H^1(\kappa,\scr{P}_{\kappa})\cong H^1(\kappa,\scr{P}_{\text{red},\kappa})$. By definition of $\scr{R}$ followed by loc.~cit., we therefore deduce that \bud H^1(\ca{O}_K, \scr{R})\cong H^1(\ca{O}_K,\scr{P}).\eud The final claim is then a consequence of \eqref{diag:claim_trivial_kernel} and the triviality over $\ca{O}_K$ of generically trivial torsors under the reductive $\ca{O}_K$-group scheme $\scr{R}$.
	\es
	\bc\label{cor:general_case_gs}
		Let $K$ be a henselian discretely valued field with a perfect residue field of characteristic $p\neq 2,3,5$ and let $\scr{P}$ be a parahoric $\ca{O}_K$-group scheme whose generic fibre $\scr{P}_K$ is a simply-connected semisimple $K$-group scheme such that $p\nmid (n+1)$ whenever $\scr{P}_K$ contains an almost simple factor of type $A_n$, for some $n\geqslant 2$.
Then, any generically trivial $\scr{P}$-torsor over $\ca{O}$ is trivial.
	\ec
	\bs
		We note that the functor $\scr{B}(-,\breve{K})$ respects product of $K$-group schemes. In particular, any parahoric model of a product of $K$-group schemes can be expressed as a product of certain parahoric models of the factors. However, by \ccite{conrad_reductive_group_schemes}{Proposition~6.4.4 and Remark~6.4.5}, the generic fibre $\scr{P}_K$ is a product of almost simple semisimple $K$-group schemes. Consequently, without loss of generality, we may assume that $\scr{P}_K$ itself is almost simple.
		\par Since $\scr{P}_K$ has no bad simple factor, the claim then follows from Theorem \ref{thm:gs_parahoric} and Corollary \ref{cor:tame_reductification_parahoric}. 
	\es
	\br\label{rem:delta>1}
	A case-by-case analysis of almost simple group schemes of type $A_n$ over henselian discretely valued fields $K$ whose residue field $\kappa$ satisfies $\cd(\kappa)\leqslant 1$, carried out in \ccite{kaletha_prasad}{\textsection10.7} (see also \cite{tits_reductive_group_local_fields} and \ccite{gross_parahoric_lecture_notes}{Proposition~5.1 and \textsection7}), shows that such a group $G$ is a bad simple factor in the sense of Definition~\ref{defn:Serre_number} if and only if $G\cong\SL_1(D)$ or $\PGL_1(D)$ for some central simple $K$-algebra $D$. Without the hypothesis $\cd(\kappa)\leqslant 1$, this classification is not known for a general perfect residue field. If such a classification were available, Proposition~\ref{prop:A_n-case} would extend Corollary~\ref{cor:general_case_gs} to the case $p\mid (n+1)$ with $p\neq 2,3,5$, even when $G$ contains an almost simple factor of type $A_n$, for some $n\geqslant 2$. 
	\er
	\br
		Our proof of Corollary~\ref{cor:general_case_gs} differs from that of Zidani \cite{zidani_gs, zidani_gs_parahoric_henselian_dvr_case}. In particular, even in the cases treated by Zidani—namely when $\scr{P}_K$ is either simply-connected, or quasi-split and adjoint—our argument provides an alternative proof under more restrictive hypotheses.
	\er

\subsection*{Acknowledgements} This project began as a conversation with Vikraman Balaji and Elden Elmanto at the Fields Institute in Toronto, during the author’s postdoctoral fellowship at the University of Toronto. The author gratefully acknowledges Elden Elmanto and Kęstutis Česnavičius for their constant support and guidance. Moreover, he owes an intellectual debt to Vikraman Balaji, Kęstutis Česnavičius, Elden Elmanto, as well as Alex Youcis, whose ideas, insights, and suggestions were central to the development of this work; many of the proofs presented here grew out of discussions with them. Special thanks are due to Philippe Gille for drawing attention to the work of Anis Zidani \cite{zidani_gs}, and to him, Ofer Gabber, Gopal Prasad, Timo Richarz, and Anis Zidani for many helpful discussions. The author also thanks Philippe Gille and Anis Zidani for a careful reading of an earlier draft and for detailed comments.
	\par This project received partial funding from PNRR grant CF 44/14.11.2022 ``\textit{Cohomological
		Hall algebras of smooth surfaces and applications}'' and NSERC Discovery grant RGPIN-2025-07114, “\textit{Motivic cohomology: theory and applications}”, and hospitality from Chennai Mathematical Institute as well as the Lodha Mathematical Sciences Institute.

\printbibliography{\let\thefootnote\relax\footnote{Simion Stoilow Institute of Mathematics of the Romanian Academy (IMAR), Bucharest, Romania\newline\hspace*{0.48cm} Email: akundu.math@gmail.com.}}{\let\thefootnote\relax\footnote{March 2025}}

\end{document}

%% file: preamble.tex
\usepackage[tmargin=2cm,bmargin=2.5cm,lmargin=1.25cm,rmargin=1.25cm]{geometry}
\usepackage[utf8]{inputenc}
\usepackage[english]{babel}
\usepackage[rm]{titlesec}
\titlelabel{\thetitle.\quad}
\usepackage[T1]{fontenc} 
\usepackage{amssymb,amsfonts}
\usepackage{amsthm}
\usepackage{tikz-cd}
\usepackage[colorinlistoftodos]{todonotes}
\usepackage[framemethod=TikZ]{mdframed} 
\usepackage{calligra,mathrsfs}
\usepackage{amsmath}
\usepackage{mathtools}
\usetikzlibrary{shapes.geometric, arrows}
\usepackage[backend=biber,style=alphabetic,sorting=nyt, maxnames=10, maxalphanames=10]{biblatex}
\usepackage{csquotes}
\usepackage{colonequals}
\usepackage{adjustbox} 
\usepackage[colorlinks, final]{hyperref} 
\usepackage{cleveref}  
\usepackage{enumerate} 
\definecolor{imperialred}{RGB}{237, 41, 57}
\definecolor{royalblue}{RGB}{64, 106, 212}
\definecolor{link}{RGB}{11,0,128}
\hypersetup{
	colorlinks = true,
	breaklinks = true,
	citecolor=royalblue,
	linkcolor=imperialred,
	urlcolor=link,
	linktocpage = true
}
\makeatletter

\makeatother

\newrobustcmd*{\citeauthorfullname}{\AtNextCite{\DeclareNameAlias{labelname}{given-family}}\citeauthor}

\newtheorem{thm}{Theorem}[section]
\newtheorem{cor}[thm]{Corollary}
\newtheorem{prop}[thm]{Proposition}
\newtheorem{lem}[thm]{Lemma}
\newtheorem{conj}[thm]{Conjecture}
\newtheorem{quest}[thm]{Question}
\newtheorem{quest*}{Question}

\theoremstyle{definition}
\newtheorem{defn}[thm]{Definition}

\newtheorem{exmp}[thm]{Example}

\theoremstyle{remark}

\newtheorem{rem}[thm]{Remark}

\newtheoremstyle{montobbo}
{7pt}
{5pt}
{}
{}%
{\bfseries}
{}%
{.5em}
{\thmnumber{\@{#1}{}~\@{#2}.}%
	\thmnote{~{\bfseries#3.}}}

\theoremstyle{montobbo}
	\newtheorem{montobo}[thm]{}
	\newcommand{\bmo}{\begin{montobo}}
	\newcommand{\emo}{\end{montobo}}
	\newtheorem{talk}{Lecture}
	\newcommand{\btalk}{\begin{talk}}
		\newcommand{\etalk}{\end{talk}}
\usepackage{footnote}
\newcommand{\ccite}[2]{\cite[#2]{#1}}
\newcommand{\stacks}[1]{\ccite{stacks-project}{\href{https://stacks.math.columbia.edu/tag/#1}{Tag \MakeUppercase{#1}}}}	

\newenvironment{manualtheorem}[1]{
	\renewcommand\thethm{#1}
	\thm
}{\endthm}
\newcommand{\but}[1]{\begin{manualtheorem}{#1}}
\newcommand{\eut}{\end{manualtheorem}}
\newenvironment{manualproposition}[1]{
	\renewcommand\thethm{#1}
	\prop
}{\endprop}
\newcommand{\bup}[1]{\begin{manualproposition}{#1}}
	\newcommand{\eup}{\end{manualproposition}}
\numberwithin{equation}{thm}
\newenvironment{teqn}{\begin{equation}\begin{tikzcd}\displaystyle}{\end{tikzcd}\end{equation}}
\newenvironment{teqn*}{\begin{equation*}\begin{tikzcd}\displaystyle}{\end{tikzcd}\end{equation*}}

\newenvironment{tpic*}{\begin{equation*}\begin{tikzpicture}}{\end{tikzpicture}\end{equation*}}

\DeclareMathAlphabet\matheu{U}{eus}{m}{n}
\DeclareMathAlphabet{\mathmm}{U}{mmambb}{m}{n}
\newcommand{\bb}[1]{\mathbb{#1}}
\newcommand{\fr}[1]{\mathfrak{#1}}
\newcommand{\ca}[1]{\mathcal{#1}}
\newcommand{\scr}[1]{\mathscr{#1}}

\newcommand{\comment}[1]{}
\newcommand{\bun}{\begin{itemize}}
\newcommand{\bn}{\begin{enumerate}}
\newcommand{\bt}{\begin{thm}}
\newcommand{\bl}{\begin{lem}}
\newcommand{\bp}{\begin{prop}}
\newcommand{\bc}{\begin{cor}}
\newcommand{\bd}{\begin{teqn}}
\newcommand{\bud}{\begin{teqn*}}	
\newcommand{\bs}{\begin{proof}}
\newcommand{\br}{\begin{rem}}
\newcommand{\bdf}{\begin{defn}}
\newcommand{\bcj}{\begin{conj}}

\newcommand{\eun}{\end{itemize}}
\newcommand{\en}{\end{enumerate}}
\newcommand{\et}{\end{thm}}
\newcommand{\el}{\end{lem}}
\newcommand{\ep}{\end{prop}}
\newcommand{\ec}{\end{cor}}
\newcommand{\ed}{\end{teqn}}
\newcommand{\eud}{\end{teqn*}}
\newcommand{\es}{\end{proof}}
\newcommand{\er}{\end{rem}}
\newcommand{\edf}{\end{defn}}
\newcommand{\ecj}{\end{conj}}
\DeclareMathOperator{\etale}{\acute{e}t}



\DeclareMathOperator{\Hom}{Hom}

\DeclareMathOperator{\Aut}{Aut}

\DeclareMathOperator{\GL}{GL}
\DeclareMathOperator{\PGL}{PGL}
\DeclareMathOperator{\SL}{SL}


\DeclareMathOperator{\Gal}{Gal}

\DeclareMathOperator{\Bl}{Bl}
\DeclareMathOperator{\aff}{aff}
\DeclareMathOperator{\spec}{Spec}



\DeclareMathOperator{\Res}{Res}
\DeclareMathOperator{\sm}{sm}

\DeclareMathOperator{\sep}{sep}
\DeclareMathOperator{\der}{der}
\DeclareMathOperator{\SC}{sc}
\DeclareMathOperator{\cd}{cd}
\DeclareMathOperator{\Stab}{Stab}
\DeclareMathOperator{\sh}{sh}


\newcommand{\tir}{\arrow[r]}

\newcommand{\iso}{\xrightarrow{\sim}}


\makeatletter\makeatother